\documentclass{gtpart}     % Basic GT/GTM/AGT style
\usepackage{pinlabel}  %%% the recommended graphics+labelling package
\usepackage{graphicx}  %%% the recommended graphics package

\title[Hyperbolic structures on link complements]
    {An alternative approach to hyperbolic structures on link complements}

\author{Morwen Thistlethwaite}
\givenname{Morwen}
\surname{Thistlethwaite}
\address{Department of Mathematics, University of Tennessee, Knoxville, TN 37996}
\email{morwen@math.utk.edu}
\urladdr{www.math.utk.edu/~morwen}

\author{Anastasiia Tsvietkova}
\givenname{Anastasiia}
\surname{Tsvietkova}
\address{Department of Mathematics, University of California, Davis, CA 95616}
\email{n.tsvet@gmail.com}
\urladdr{www.math.ucdavis.edu/~tsvietkova}

\keyword{Classical Link}
\keyword{Hyperbolic Structure}
\keyword{Tangle}
\subject{primary}{msc2000}{57M25}
\subject{secondary}{msc2000}{57M50}

\volumenumber{}
\issuenumber{}
\publicationyear{}
\papernumber{}
\startpage{}
\endpage{}
\doi{}
\MR{}
\Zbl{}
\received{}
\revised{}
\accepted{}
\published{}
\publishedonline{}
\proposed{}
\seconded{}
\corresponding{}
\editor{}
\version{}

\newtheorem{thm}{Theorem}[section]    % Standard theorem environment
\newtheorem{prop}[thm]{Proposition}
\newtheorem{conj}[thm]{Conjecture}
\newtheorem{lemma}{Lemma}[thm]
\newtheorem{cor}{Corollary}[thm]
\theoremstyle{definition}
\newtheorem{defn}[thm]{Definition}    % Definition environment with 
\newtheorem{note}[thm]{Note}
\newtheorem*{rem}{Remark}             % Unnumbered environment for remarks.
\newtheorem{example}{Example}[thm]
\def\bc{\,,\,}
\def\h{\hspace{.5em}}
\def\ha{\hspace{.3em}}
\def\phn{\phantom{--}}
\def\phl{\phantom{l}}
\def\dis{\displaystyle}
\def\v{\vspace{.5ex}}
\def\vv{\vspace{1ex}}

\begin{document}

\begin{abstract}    % type your abstract below

An alternative method is described
for determining the hyperbolic structure on a link complement, and some of its elementary
consequences are examined.  The method is particularly suited to alternating links.

\end{abstract}

\maketitle

%%%%%%%%%%%%%%%%%%%%   Start of main body of article

\begin{section}{Overview}

The purpose of this article is to describe an alternative method for calculating the
hyperbolic structure on a classical link complement.  The method does not use
an ideal triangulation of the complement, but instead considers the shapes of
ideal polygons bounding the regions of a diagram of the link.  In order to
guarantee the applicability of our method, we shall impose a
``minimality'' condition on the checkerboard surfaces of our link diagrams:

\begin{defn}
A diagram of a hyperbolic link is {\it taut} if
each associated checkerboard surface is incompressible and boundary incompressible
in the link complement, and moreover does not contain any simple closed curve
representing an accidental parabolic.
\end{defn}

From this definition it follows that if \( \alpha \) is a proper, non-separating
arc in a checkerboard surface associated to a taut diagram, and \( \widetilde{\alpha} \)
is a lift of \( \alpha \) to the universal cover \( \mathbb H^3 \), then the ends
of \( \widetilde{\alpha} \) are at the centres of distinct horoballs; thus \( \alpha \) is properly
homotopic to a geodesic.  In particular, at each crossing of the diagram, the arc
travelling vertically from underpass to overpass, {\it i.e.} a ``polar axis'' in
the terminology of \cite{Men}, gives rise in this manner to a geodesic; such geodesics,
henceforth called {\it crossing geodesics}, will form the edges of the ideal polygons
mentioned above.

Although the method is applicable to any taut link diagram, we are particularly interested
in applying it to hyperbolic alternating links, as the resulting hyperbolicity equations
assume a reasonably pleasing form.  We recall that it is proved in \cite{Men} that
prime alternating link complements cannot contain
essential tori, and since the only alternating torus links are those of type $(2,n)$\,,
it follows from W. Thurston's hyperbolization theorem that an alternating link is
hyperbolic if and only if it is prime and is not a $(2,n)$--torus link.  From \cite{Men} a
reduced alternating link diagram represents a prime link if and only if it is prime in
the diagrammatic sense, and from \cite{MT}  each reduced alternating diagram of a
$(2,n)$--torus link is standard; therefore one can tell by inspection whether a link
presented as a reduced alternating diagram is hyperbolic.

\begin{prop}  Each reduced alternating diagram of a hyperbolic
alternating link is taut.
\end{prop}

\begin{proof}
It is proved in \cite{MT} that the checkerboard surfaces for
such link diagrams are incompressible and boundary incompressible,
and it is proved in \cite{Ad2,FKP} that they are quasi-fuchsian,
hence contain no accidental parabolics.
\end{proof}

\begin{note}
In \cite{FKP}, the authors state their results for a more
general class of diagrams than alternating.  However, the spanning surfaces
considered are so-called {\it state surfaces}; for non-alternating diagrams
these are different from checkerboard surfaces, and we do not know at present whether
these can be incorporated into our method of computing hyperbolic structures.
\end{note}

It would be interesting to know whether there exists a hyperbolic link not admitting
a taut diagram.
\end{section}

\begin{section}{The geometry of an ideal polygon}

Let \( F \) be a checkerboard surface for a connected diagram \( D \) of a link \( L \).
Then \( F \) is the union of disks, one for each region coloured say black in the checkerboard
colouring of the diagram.  The boundary of each disk is an alternating sequence of (i) sub-arcs of
the link travelling between adjacent crossings incident to the region, and (ii)
``polar axis'' arcs travelling between the underpass and the overpass at a crossing.
The disks are glued together along the polar axis arcs.

Now suppose that \( D \) is taut; let \( R \) be a black region of \( D \) with \( n \geq 2 \) sides,
and let \( \Delta_R \subset F \) be the associated disk.  Then \( \Delta_R - L \) is homeomorphic to a
disk with \( n \) points of its boundary removed, which we may describe as a ``filled-in ideal
$n$--gon''; this lifts homeomorphically to a filled-in ideal $n$--gon
\( \widetilde{\Delta_R} \) in the upper half-space model of \( \mathbb H^3 \).
The \( n \) ideal vertices of \( \widetilde{\Delta_R} \) correspond to the \( n \) arcs
of \( \Delta_R \cap L \) (which in turn correspond to edges of the region \( R \)), and the edges
of the ideal $n$--gon boundary of \( \widetilde{\Delta_R} \)
are lifts of the interiors of the \( n \) polar axis arcs in \( \partial \Delta_R \).

In order to proceed further, we need to show that \( \widetilde{\Delta_R} \) satisfies a
non-degeneracy condition.

\begin{prop}  The \( n \) ideal vertices of \( \widetilde{\Delta_R} \) are
pairwise distinct.
\end{prop}

\begin{proof}
Let \( \alpha_1 \bc \alpha_2 \) be any two arc components of \( \Delta_R \cap L \),
and let \( \gamma \) be an arc properly embedded in \( \Delta_R \) that travels from a point of
\( \alpha_1 \) to a point of \( \alpha_2 \).  Then the interior of \( \gamma \) lifts to an
arc in \( \widetilde{\Delta_R} \) travelling between the corresponding ideal vertices.
Since the link diagram is taut and \( \gamma \) is non-separating in the checkerboard
surface \( F \), the conclusion follows.
\end{proof}

As usual, we identify the boundary of \( \mathbb H^3 \) with the Riemann sphere
\( \mathbb C \cup \{ \infty \} \).  Let \( R \) be a region of the link diagram with
at least three sides, and let the ideal vertices of \( \widetilde{\Delta_R} \)
be \( z_1 \bc \dots \bc z_n \) in cyclic order; then, from Proposition 2.1,
these \( n \) points define an ideal $n$--gon \( \widetilde{\Pi_R} \) in \( \mathbb H^3 \),
with geodesic edges that are pairwise distinct.

Let \( \gamma_i \) be the geodesic edge of \( \widetilde{\Pi_R} \) joining \( z_i \) with
\( z_{i+1} \) (where indices are taken modulo $n$).  We define the {\it shape parameter}
\( \zeta_i \) of \( \gamma_i \) to be the cross-ratio
\[ \zeta_i \;=\; \frac{(z_{i-1} - z_i)(z_{i+1} - z_{i+2})}{(z_{i-1} - z_{i+1})(z_i - z_{i+2})} \h,\]
with the usual rules about cancelling \( \pm \infty \) terms.  If we perform an isometry
of \( \mathbb H^3 \) to place the vertices \( z_{i-1} \bc z_i \bc z_{i+1} \) at \( 1 \bc \infty \bc 0 \)
respectively, then the vertex \( z_{i+2} \) will be placed at \( \zeta_i \), and
we see that the collection of $n$ shape parameters determines the isometry class
of the ideal $n$--gon.

It follows that for a 3--sided polygon each shape parameter is equal to \( 1 \);
it is also easy to check that for a 4--sided polygon the sum of two consecutive
shape parameters is \( 1 \), whence opposite shape parameters are equal.  For general $n$,
we may obtain convenient equations relating the \( \zeta_i \) from the fact that the polygon
closes up.  Specifically, if we place the polygon so that
\( z_{i-1} = 1 \bc z_i = \infty \bc z_{i+1} = 0 \), then the isometry \( \psi_i \) given by the
M\"{o}bius transformation \(\dis z \mapsto \frac{-\zeta_i}{z - 1} \) maps
\( z_{i-1} \bc z_i \bc z_{i+1} \) to \( z_i \bc z_{i+1} \bc z_{i+2} \) respectively.
Since the polygon \( \widetilde{\Pi_R} \) closes up, the composite
\( \psi_n \circ \dots \circ \psi_2 \circ \psi_1 \) must equal the identity, and passing to matrices,
we see that we have an identity
\begin{equation}
   \left[\begin{array}{cc}0&-\zeta_n\\1&-1\end{array}\right] \h\dots\h
   \left[\begin{array}{cc}0&-\zeta_2\\1&-1\end{array}\right]
   \left[\begin{array}{cc}0&-\zeta_1\\1&-1\end{array}\right] \quad\sim\quad
   \left[\begin{array}{cc}1&0\\0&1\end{array}\right] \h,
\end{equation}
where \( \sim \) denotes equality up to multiplication by a non-zero scalar matrix.
From the \( (2 \bc 1)$--entry of this product
we can read off a polynomial relation \( f_n = 0 \) in the \( \zeta_i \).  It is then
easily checked, using induction on \( n \), that the polynomials \( f_n \) may be defined
recursively by
\begin{equation}
f_3 \equiv 1 - \zeta_2 \quad,\quad f_4 \equiv 1 - \zeta_2 - \zeta_3 \quad,\quad
f_n \equiv f_{n-1} - \zeta_n f_{n-2} \h (n \geq 5) \quad.
\end{equation}
We observe that the polynomial \( f_n \) involves the \( n - 2 \)
shape parameters \( \zeta_2 \bc \zeta_3 \bc \dots \bc \zeta_{n-1} \), and that
\( f_n \) is of degree 1 in each of these shape parameters.  In particular,
the relation \( f_n = 0 \) allows one to express each of these \( n - 2 \) shape
parameters as a function of the other \( n - 3 \).  Let
\( f_n^+ \h (f_n^-) \) be the polynomial obtained from \( f_n \) by
increasing all indices by 1 (resp. decreasing all indices by 1).
Then \( f_n^+ \) is independent from \( f_n \), as it is the only one of the two
that involves \( \zeta_n \);
also, \( f_n^- \) is independent from both \( f_n \) and \( f_n^+ \), as it is the
only one of the three that involves \( \zeta_1 \).
In fact
\( \{ f_n^- = 0 \bc f_n = 0 \bc f_n^+ = 0 \} \) must be a complete set of relations
for the $n$ shape parameters of a generic ideal $n$--gon, as the triple transitivity
of the action of the group of M\"{o}bius transformations on the boundary \( \mathbb C \cup \{ \infty \} \)
dictates that the isometry class of an ideal polygon with \( n \)
sides has \( n - 3 \) geometric degrees of freedom.

It is immediate from the definition of shape parameter that an ideal polygon lies in a
hyperbolic plane if and only if all its shape parameters are real.
For highly symmetric links, ideal polygons are often encountered that are {\it regular}, in the
sense that all \( \zeta_i \) are equal.

\begin{note}
In order to conform to various sign conventions, later we shall
be obliged to deal with the {\it complex conjugates} of shape parameters.  However, since
complex conjugation is a field automorphism of the complex numbers, these complex conjugates
\( \overline{\zeta_i} \) satisfy the same polynomial relations as the \( \zeta_i \).
\end{note}

\begin{prop}
The common shape parameter for a regular $n$--sided ideal polygon is\\
\( \frac{1}{4}\sec^2\frac{\pi}{n} \).
\end{prop}

\begin{proof}
We may assume that the $n$ ideal vertices of the polygon are evenly spaced
around a unit circle in \( \mathbb C \); specifically, we assume that the ideal vertices are
\( w^i \h (0 \leq i \leq n-1) \), where \( w = e^{2\pi i/n} \).  The cross-ratio of the first four
of these points is then
\[ \frac{(1 - w)(w^2 - w^3)}{(1 - w^2)(w - w^3)} \h=\h \frac{w^2(1 - w)^2}{w(1-w^2)^2} \h=\h \frac{w}{(1 + w)^2} \h.\]
Noting that the line segment joining \( 0 \) with \( 1 + w \) is a diagonal of the rhombus with vertices
\( 0 \bc 1 \bc 1 + w \bc w \), we see that the modulus of \( 1 + w \) is \(\dis 2 \cos \pi/n \).  Also, since
this diagonal bisects the angle of the rhombus at \( 0 \), the argument of the above cross-ratio is \( 0 \);
the result follows.
\end{proof}

It follows from Proposition 2.3 that for regular polygons \( \zeta \) decreases monotonically to the limit
\(\frac{1}{4} \) as \( n \to \infty \).

In the absence of symmetry, polygons with four or more sides need not be regular, and
need not lie in a hyperbolic plane, although for alternating links it seems from experiment
that they are close to being planar, and never deviate very far from being regular.  Here are
two examples of this phenomenon.

\begin{figure}[ht!]
\labellist
\hair 2pt
\pinlabel $9a37$ at 195 -30
\pinlabel $11a79$ at 705 -30
\endlabellist
\centering

\includegraphics[scale=.4]{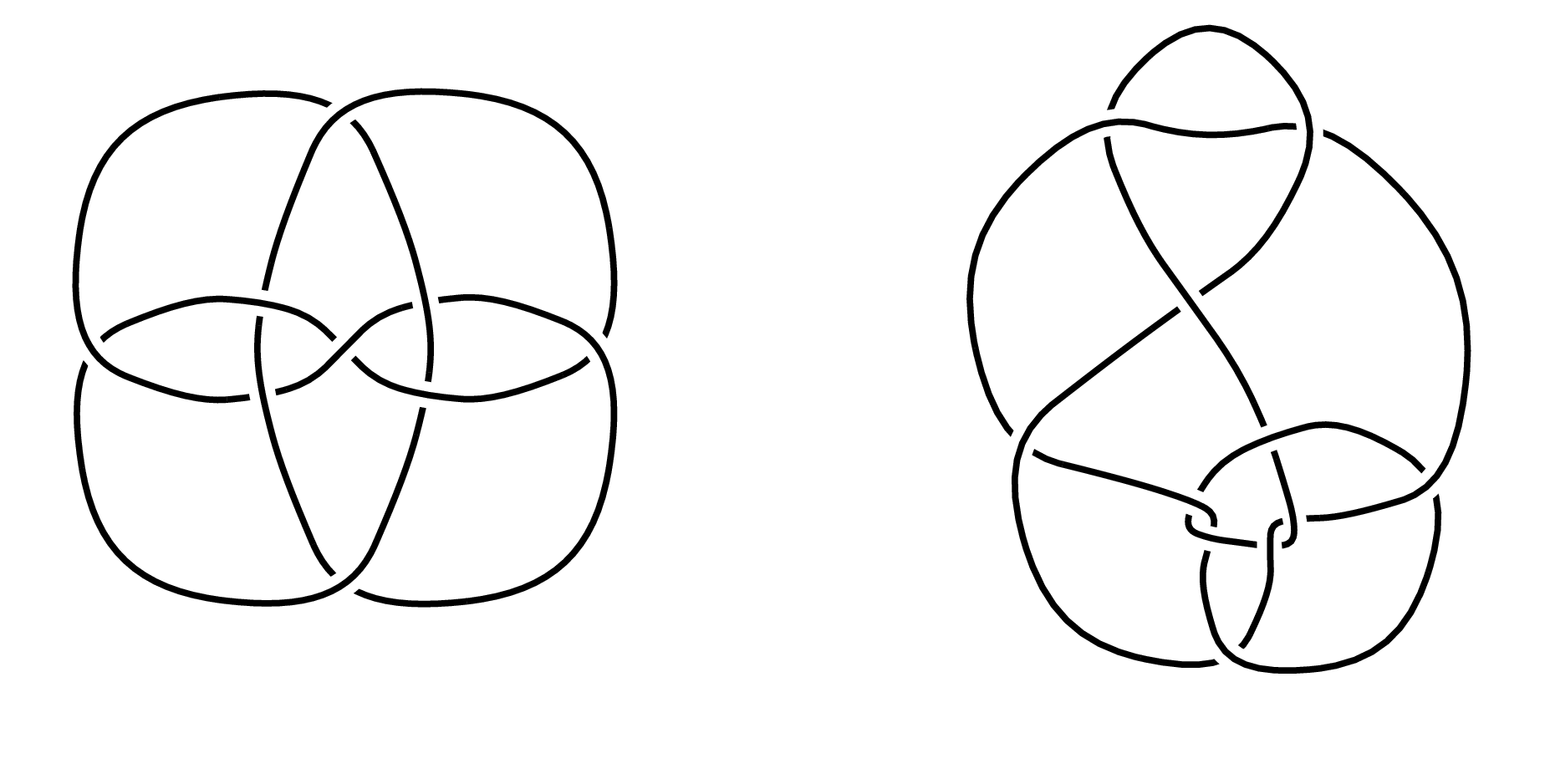}
\caption{}
\end{figure}

\begin{example}
{\it The knot \( 9a37 \) in the Dowker-Thistlethwaite listing.}
\end{example}

There is a symmetry of order \( 3 \) cyclically permuting the three regions with four sides.
For each of these regions, the shape parameters corresponding to the ``north'' and ``south''
crossings are both \( 0.469789 - 0.090643\,i \), and the shape parameters corresponding
to the ``east'' and ``west'' crossings are both \( 0.530211 + 0.090643\,i \), these values
being rounded to six decimal places.  The imaginary parts are seen to be quite small, and the
shape parameters are fairly close to that of a regular 4--sided region, namely \( \frac{1}{2} \).

\begin{example}
{\it The knot \( 11a79 \) in the Dowker-Thistlethwaite listing.}
\end{example}

For the small 5--sided region in the lower-middle part of the diagram, we begin at the top crossing
and proceed around the region in a counterclockwise direction.  To six decimal places, the five shape parameters
for this region are as follows:
\h  \( 0.312331 - 0.008243\,i \bc\, 0.449632 - 0.007097\,i \bc\, 0.346369 + 0.018155\,i \bc\,
0.370339 - 0.024868\,i \bc\, 0.432793 + 0.022291\,i \).  This time we compare with the shape parameter of
a regular 5--sided ideal polygon, \( \frac{1}{4}\sec^2\frac{\pi}{5} = (3 - \sqrt{5})/2 \approx 0.381966 \).

In the next section we shall see how to use the peripheral structure of the link
complement to set up a system of equations for determining the shape parameters of
the ideal polygons and for determining how the polygons are situated relative to one another.
The unknowns of these equations will be complex numbers attached to the edges and
crossings of the diagram; these complex number ``labels'' will in fact determine the
complete hyperbolic structure of the link complement.
\end{section}

\begin{section}{Edge and crossing labels}

We assume throughout that horospherical cross-sections of the cusps have been chosen
so that a (geodesic) meridian curve on the cross-sectional torus has length \( 1 \).
This guarantees \cite{Ad1} that cross-sectional tori from distinct cusps are disjoint,
and that each torus is
embedded in the link complement, with the exception of the figure-eight knot complement,
where the cross-sectional torus touches itself in two points.

The preimage of each cross-sectional torus in the universal cover \( \mathbb H^3 \)
is a union of horospheres, and we specify a complex affine structure on each horosphere by
declaring (for convenience) that meridional translation is through unit distance in
the positive real direction.  We also would like the translation corresponding to a
longitude on the torus to have positive imaginary part, so in order to keep the standard
orientation of the complex plane and the usual ``right-hand screw'' convention relating
the directions of meridian and longitude, we view the torus {\it from the thick part}
of the manifold, {\it i.e.} from the opposite side to the cusp.

We assume that coordinates are chosen so that one of
the horospheres is the Euclidean plane \( H_\infty \) of (Euclidean) height \( 1 \) above
the $xy$--plane, and that it has the standard affine structure; thus on that horosphere
the meridional translation is represented by the matrix
\( \left[\begin{array}{cc}1&1\\0&1\end{array}\right] \).

\begin{figure}[ht!]
\labellist
\small
\pinlabel $\mathrm{H_{\, i-1}}$ at -59 38
\pinlabel $\mathrm{z_{\, i-1}}$ at 1 8
\pinlabel $\mathrm{P_{\, i-1}}$ at 31 133
\pinlabel $\mathrm{\gamma_{\: i-1}}$ at 84 130
\pinlabel $\mathrm{H_{\, i}}$ at 95 21
\pinlabel $\mathrm{Q_{\, i}}$ at 129 99
\pinlabel $\mathrm{z_{\, i}}$ at 151 -17
\pinlabel $\mathrm{P_{\, i}}$ at 156 99
\pinlabel $\mathrm{\gamma_{\, i}}$ at 192 154
\pinlabel $\mathrm{H_{\, i+1}}$ at 247 76
\pinlabel $\mathrm{Q_{\, i+1}}$ at 257 172
\pinlabel $\mathrm{z_{\, i+1}}$ at 335 31
\endlabellist
\centering

\includegraphics[scale=0.7]{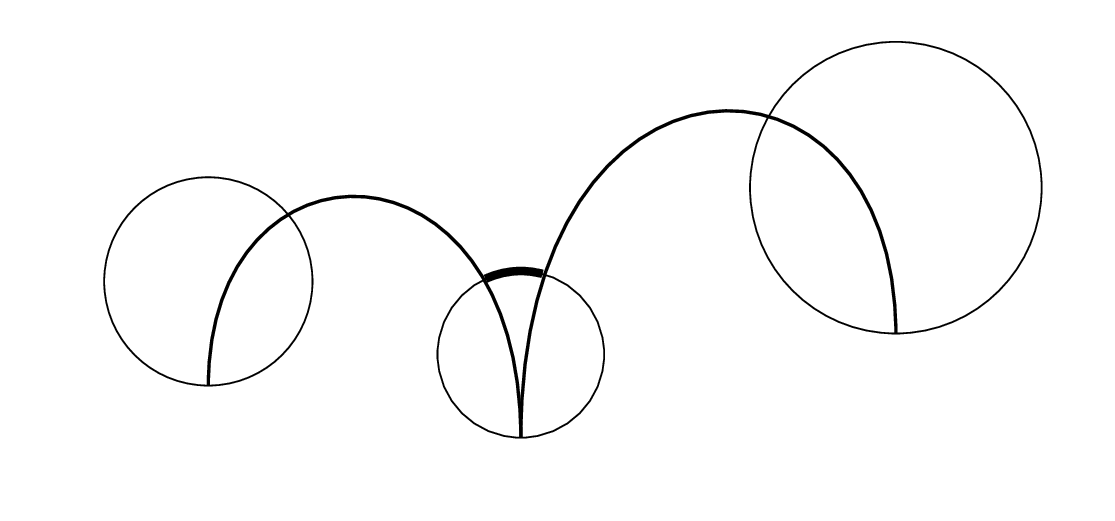}
\caption{}
\end{figure}

Now consider an ideal polygon \( \widetilde{\Pi_R} \) in \( \mathbb H^3 \), associated
to a region \( R \) of the link diagram, and
with ideal vertices \( z_1 \bc \dots \bc z_n \) as above.  Each vertex \( z_i \) is
the center of a horosphere \( H_i \), and each geodesic edge \( \gamma_i \) of
\( \widetilde{\Pi_R} \) meets \( H_i \bc H_{i+1} \) in points \( P_i \bc Q_{i+1} \)
respectively (Fig.\ha 2).  At this stage, for book-keeping purposes we need to choose an
orientation of the link, and this orientation will determine a direction on the
geodesic arc on \( H_i \) joining \( Q_i \) with \( P_i \).  From the affine
structure on \( H_i \) we now have a complex number determining a translation
mapping one of \( P_i \bc Q_i \) to the other, depending on this direction.
This complex number is affixed to the side of the corresponding edge \( E \) of the link
diagram incident to the region \( R \), and will be called an {\it edge label}.

Next we note that there is a simple relation between the edge labels on the two sides
of an edge \( E \) of the link diagram.  Let the regions incident to \( E \) be
\( R \bc S \), and let the corresponding labels affixed to the two sides of \( E \)
be \( u_R^E \bc u_S^E \) respectively.  These labels correspond to geodesic arcs on
a horosphere that descend to arcs \( \alpha_R^E \bc \alpha_S^E \) on the peripheral
torus, joining the points of intersection of the torus with two successive crossing
geodesics.  Let \( \mu \) denote a meridian curve on the torus, oriented as usual
according to the ``right-hand screw'' rule, and for an arc \( \alpha \) let
\( \overline{\vphantom{\text{A}}\alpha} \) denote its reverse.  The loop \( \alpha_R^E * \overline{\alpha_S^E} \)
is homotopic to \( \kappa\mu \), where
\[ \kappa \h=\h \begin{cases}
   \phn \phantom{-}1 \quad &\text{if \( E \) ascends from left to right}\\
   \phn {-} 1      &\text{if \( E \) descends from left to right}\\
   \phn \phantom{-}0       &\text{if \( E \) is level}
   \end{cases} \]
as one looks from the interior of the region \( R \).  It follows that
\( u_R^E - u_S^E = \kappa \), where \( \kappa \) is as above (see Fig.\ha 3, which illustrates
the case \( \kappa = 1 \)).

\begin{figure}[ht!]
\labellist
\pinlabel $\mathrm{\gamma_{\, 1}}$ at -267 -36
\pinlabel $\mathrm{\mu}$ at -202 -144
\pinlabel {\scriptsize$\mathrm{S}$} at -137 -6
\pinlabel $\mathrm{u_{\hskip.5pt R}}$ at -17 -70
\pinlabel $\mathrm{u_{\hskip.5pt S}}$ at 41 -3
\pinlabel $\mathrm{R}$ at 156 -115
\pinlabel $\mathrm{\gamma_{\, 2}}$ at 268 -63
\endlabellist
\centering

\includegraphics[scale=0.52]{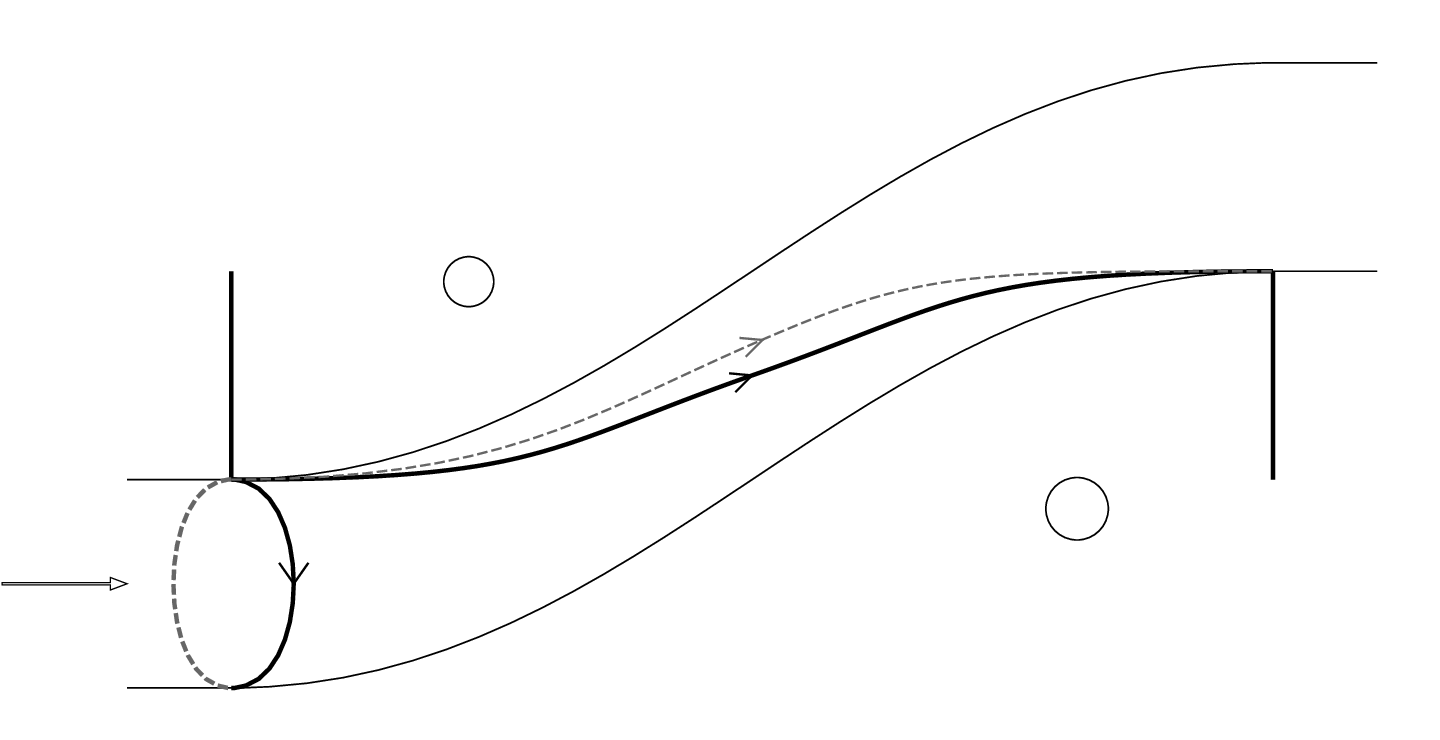}
\caption{}
\end{figure}

In the case of reduced alternating diagrams, the relation between \( u_R^E \) and \( u_S^E \)
is particularly simple.  If we colour the regions in checkerboard fashion, then the view from
inside regions of one colour has all edges on the boundary of the region ascending from left
to right, and the view from regions of the other colour has edges descending from left to right.
We may take the convention that regions of the former type are black, and those of the
latter type white.  Then, if \( R \) is a black region and \( S \) is an adjacent white
region, we have \( u_R^E - u_S^E = 1 \).

We turn now to crossing labels.  To each crossing geodesic \( \gamma \) we assign a
complex number \( w_\gamma \) as follows.  We lift \( \gamma \) to a geodesic
\( \widetilde{\gamma} \) in \( \mathbb H^3 \) joining the centers of horospheres
\( H_1 \bc H_2 \)  The meridional direction on \( H_i \h (i = 1 \bc 2) \) together
with the geodesic \( \widetilde{\gamma} \) defines a
hyperbolic half-plane \( \Sigma_i \) containing \( \widetilde{\gamma} \), and we define
the argument of \( w_\gamma \) by \( \arg w_\gamma = \phi_\gamma + \pi \),
where \( \phi_\gamma \) is the angle between the half-planes \( \Sigma_1 \bc \Sigma_2 \),
the sign of \( \phi_\gamma \) being determined by the convention that a right-handed screw is
positive.  Thus we have \( \arg(-w_\gamma) = \phi_\gamma \); the reason for choosing
\( -w_\gamma \) here rather than \( w_\gamma \) is explained in the next paragraph.
In essence we are defining the angle \( \phi_\gamma \) as the angle between the
two meridional directions by parallel transport along the geodesic \( \widetilde{\gamma} \).
The specification of \( w_\gamma \) is completed by defining its modulus to be \( e^{-d} \),
where \( d \) is the hyperbolic length of the part of \( \widetilde{\gamma} \) between
\( H_1 \) and \( H_2 \).

If we take \( H_1 \) to be the horosphere \( H_\infty \) defined above, then \( |w_\gamma| \)
is the Euclidean diameter of the horosphere \( H_2 \).  The argument of \( w_\gamma \) may
loosely be interpreted as the angle (in the hyperbolic structure) between the overpass and
underpass at the crossing in question.  The reason for choosing \( \phi_\gamma + \pi \)
rather than \( \phi_\gamma \) as the argument of \( w_\gamma \) is illustrated in Fig.\ha 4.
In that figure are illustrated an overpass and underpass that are parallel; thus, if
\( \arg(w_\gamma) \) is to measure the angle between these strands, we would like
\( \arg(w_\gamma) \) to be zero.  However, at the top end of the intercusp segment of
the geodesic, the direction of the meridian is away from the viewer, whereas at the
bottom end the direction of the meridian is towards the viewer.

\begin{figure}[ht!]
\centering

\includegraphics[scale=0.85]{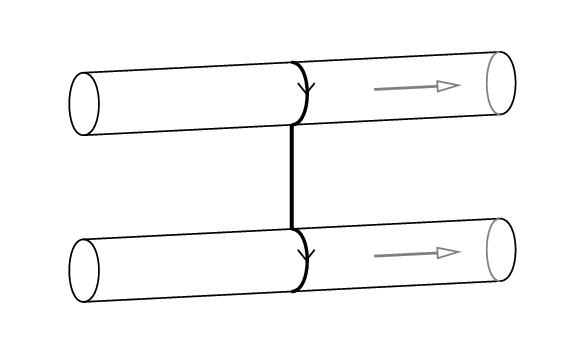}
\caption{}
\end{figure}

We note that the isometry represented by the
matrix \( \left[\begin{array}{cc}0&-w_\gamma\\1&0\end{array}\right] \) maps \( H_1 \) to
\( H_2 \), respecting affine structures.

In the degenerate case of a two-sided region, the two arcs travelling from overpass to
underpass of the two crossings are parallel, hence are homotopic to the same geodesic.
Therefore the two crossing labels are equal, the ideal polygon for this region
collapses, and the two edge labels for the region are \( 0 \).  For taut
diagrams, in particular reduced alternating diagrams, we have the following converse:
\begin{thm}
Let \( R \) be a region of at least three sides of a taut
diagram \( D \) of a hyperbolic link \( L \).  Then each edge label for \( R \) is non-zero.
\end{thm}

\begin{proof}
Suppose that \( E \) is some edge of the region \( R \) for which the
label \( u_R^E \) is zero.  Let \( z \) be the ideal vertex of \( \widetilde{\Pi_R} \)
corresponding to this edge, and let \( H \) be the horosphere centred at \( z \).
Since \( u_R^E = 0 \), the two geodesic edges of \( \widetilde{\Pi_R} \) issuing from \( z \)
meet \( H \) in the same point, hence are equal, contradicting Proposition 2.1.
\end{proof}

\begin{conj}
Let \( D \) be a reduced alternating diagram of a hyperbolic link.
Each edge label for \( D \) has non-negative imaginary part.
\end{conj}

The above definition of crossing labels applies to any geodesic \( \gamma \) travelling from
cusp to cusp.  Indeed, such a geodesic may be regarded as belonging to a crossing of
some diagram of the link.

We conclude this section with a remark on symmetries.  Let \( h: (S^3 \bc L) \to (S^3 \bc L) \)
be a homeomorphism; then, modifying \( h \) by a homotopy if necessary, we may assume
that the restriction of \( h \) to \( S^3 - L \) is  an isometry.  Let us suppose
that \( h \) maps a crossing geodesic \( \gamma \) to \( \gamma' \).  Then the moduli
of \( w_\gamma \bc w_{\gamma'} \) will be equal.  If \( h \) preserves the orientation of \( S^3 \),
the associated labels \( w_\gamma \bc w_{\gamma'} \) will also have equal arguments,
whence \( w_{\gamma'} = w_{\gamma} \); on the other hand,
if \( h \) reverses the orientation of \( S^3 - L \) the argument of
\( w_\gamma \) will be negated, whence \( w_{\gamma'} = \overline{w_\gamma} \).
Edge labels are affected similarly under the action of \( h \):  if \( h \) maps a diagram of \( L \)
to itself and \( u \bc u' \) are edge labels that correspond under the homeomorphism, then
\( u' = u \) or \( u' = \overline{u} \) depending on whether \( h \) preserves or reverses the
orientation of \( S^3 \).
\end{section}

\begin{section}{The hyperbolicity equations in the edge and crossing labels}

We begin with a useful identity. Let \( R \) be a region of a link diagram with at least
three sides; from Theorem 3.1 this condition ensures that all edge labels for \( R \)
are non-zero.  Let \( \gamma \) be a geodesic edge of the ideal polygon bounding \( R \),
let \( w \) be the crossing label attached to \( \gamma \),
and let \( u_R^{E_1} \bc u_R^{E_2} \) be the labels for the edges incident
to this crossing, the suffix \( R \) indicating of course that the labels are
placed on the sides of these edges in the region \( R \).  Then we have the following
relation between the shape parameter \( \zeta_\gamma \) and the edge labels
\( u_R^{E_1} \bc u_R^{E_2} \):
\begin{prop}

\begin{equation}
\overline{\zeta_\gamma} \h=\h \frac{\kappa w}{u_R^{E_1} u_R^{E_2}} \quad,
\end{equation}
where \( \kappa = 1 \) if one edge is directed towards the crossing and one away from
the crossing, and where \( \kappa = -1 \) if both edges are directed towards the
crossing or both away from the crossing.
\end{prop}

\begin{rem}
The presence of complex conjugation in \( \overline{\zeta_\gamma} \)
is merely an artefact of our various sign conventions; as observed in Note 2.2, it
is not an obstacle, as the \( \overline{\zeta_i} \) satisfy the same relations as
the \( \zeta_i \).  By means of formula {\bf (3)}, each equation in (complex conjugates of)
shape parameters may now be regarded as an equation in edge and crossing labels; moreover,
by clearing denominators, we may regard these equations as polynomial equations in the
labels.
\end{rem}

\begin{lemma}
Fig.\ha 5 illustrates a configuration in the hyperbolic plane containing
a horocycle of Euclidean diameter \( D \), and geodesics \( \gamma_1 \bc \gamma_2 \) originating
from the centre of the horocycle and meeting the horocycle
at points \( P \bc Q \) respectively.  The ends of \( \gamma_2 \) are Euclidean
distance \( d \) apart.  Then the hyperbolic distance along the horocycle from \( P \) to \( Q \)
is \( D/d \).
\end{lemma}

\begin{proof}
To obtain the hyperbolic distance along the horocycle from \( P \) to \( Q \)
we perform the isometry that is inversion in the circle \( C \) illustrated, with Euclidean centre
at the foot of \( \gamma_1 \) and with radius \( d \).  This inversion maps \( \gamma_1 \) into
itself, maps \( \gamma_2 \) to the geodesic \( \gamma_2' \) illustrated, and maps the horocycle
to the Euclidean horizontal line at height \( d^2/D \) above the boundary.  The points \( P \bc Q \)
are mapped to points \( P' \bc Q' \) respectively on this image horocycle, and the circular arc
on the horocycle joining \( P \bc Q \) is mapped to the (horizontal) Euclidean line segment joining
\( P' \bc Q' \).  This line segment has hyperbolic length \( d / (d^2/D) = D/d \), QED.
\end{proof}

\begin{figure}[ht!]
\labellist
\pinlabel $\mathrm{C}$ at -102 87
\pinlabel $\mathrm{D}$ at -60 33
\pinlabel $\mathrm{P}$ at -8 72
\pinlabel $\mathrm{P'}$ at -11 243
\pinlabel $\mathrm{\gamma_{\, 1}}$ at 15 168
\pinlabel $\mathrm{Q}$ at 25 63
\pinlabel $\mathrm{\gamma_{\, 2}}$ at 52 72
\pinlabel $\mathrm{d}$ at 58 -30
\pinlabel $\mathrm{\gamma\,'_{\, 2}}$ at 137 133
\pinlabel $\mathrm{Q\,'}$ at 135 243
\endlabellist
\centering

\includegraphics[scale=0.75]{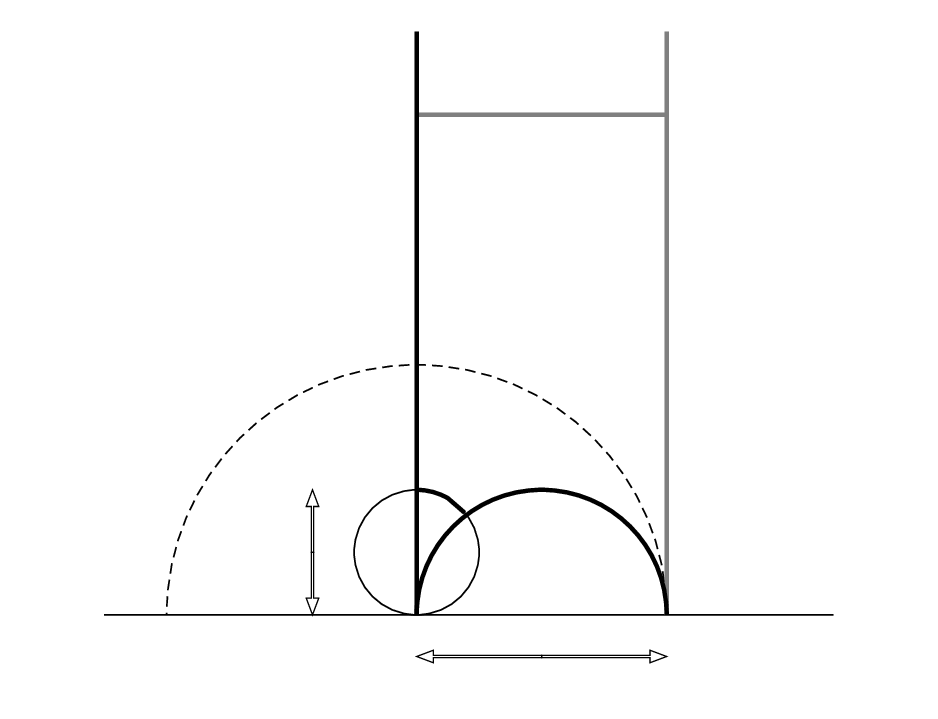}
\caption{}
\end{figure}

\begin{proof}[Proof of Proposition 4.1]

We shall verify {\bf (3)} in the case where both edges are directed away from the crossing,
in which case \( \kappa = -1 \).  The other cases are then verified simply by changing the sign
of one or both of the edge labels.

Let \( z_0 \bc z_1 \bc z_2 \bc z_3 \)  be consecutive ideal vertices
of the ideal polygon corresponding to the region \( R \), such that the geodesic
\( \gamma \) joins \( z_1 \) to \( z_2 \).  Applying a suitable isometry of hyperbolic
space, we may assume that \( z_0 = |u_R^{E_1}| \bc z_1 = \infty \bc z_2 = 0 \).  The horosphere
of which \( z_1 \) is centre is the Euclidean plane at height 1 above the boundary plane,
and the horosphere of which \( z_2 \) is centre has diameter \( |w| \).  Cancelling infinite
terms in the usual way, the cross-ratio defining \( \zeta_\gamma \) is seen to be
\( \frac{z_3}{z_0} \).  The situation is illustrated in Fig.\ha 6, where \( \mu_1 \bc \mu_2 \)
denote meridian vectors, and (for notational simplicity) \( u_1 \bc u_2 \) denote
\( u_R^{E_1} \bc u_R^{E_2} \) respectively.

First we check that the moduli of the two sides of {\bf (3)} agree.  This follows from
Lemma 4.1.1, applied to the hyperbolic plane containing \( z_1 \bc z_2 \bc z_3 \), and
where the quantities of the lemma apply as follows:
\( \gamma_1 = \gamma \;,\; D = |w| \;,\; d = |z_3| = |\zeta_\gamma| |u_R^{E_1}| \),\,
and the distance from \( P \) to \( Q \) along the horocycle is \( |u_R^{E_2}| \).

\begin{figure}[ht!]
\labellist
\pinlabel $\mathrm{\gamma_{\, 1}}$ at 19 172
\pinlabel $\mathrm{\mu_{2}}$ at 23 -2
\pinlabel $\mathrm{z_{\hskip.6pt 2} = 0}$ at 37 -75
\pinlabel $\mathrm{u_{\hskip.3pt 2}}$ at 70 26
\pinlabel $\mathrm{\mu_{1}}$ at 85 121
\pinlabel $\mathrm{\gamma_{\, 2}}$ at 124 -16
\pinlabel $\mathrm{z_{\hskip.6pt 3}}$ at 150 -112
\pinlabel $\mathrm{u_{\hskip.3pt 1}}$ at 165 77
\pinlabel $\mathrm{z_{\hskip.6pt 0}}$ at 280 -77
\pinlabel $\mathrm{\gamma_{\, 0}}$ at 299 171
\endlabellist
\centering

\includegraphics[scale=0.625]{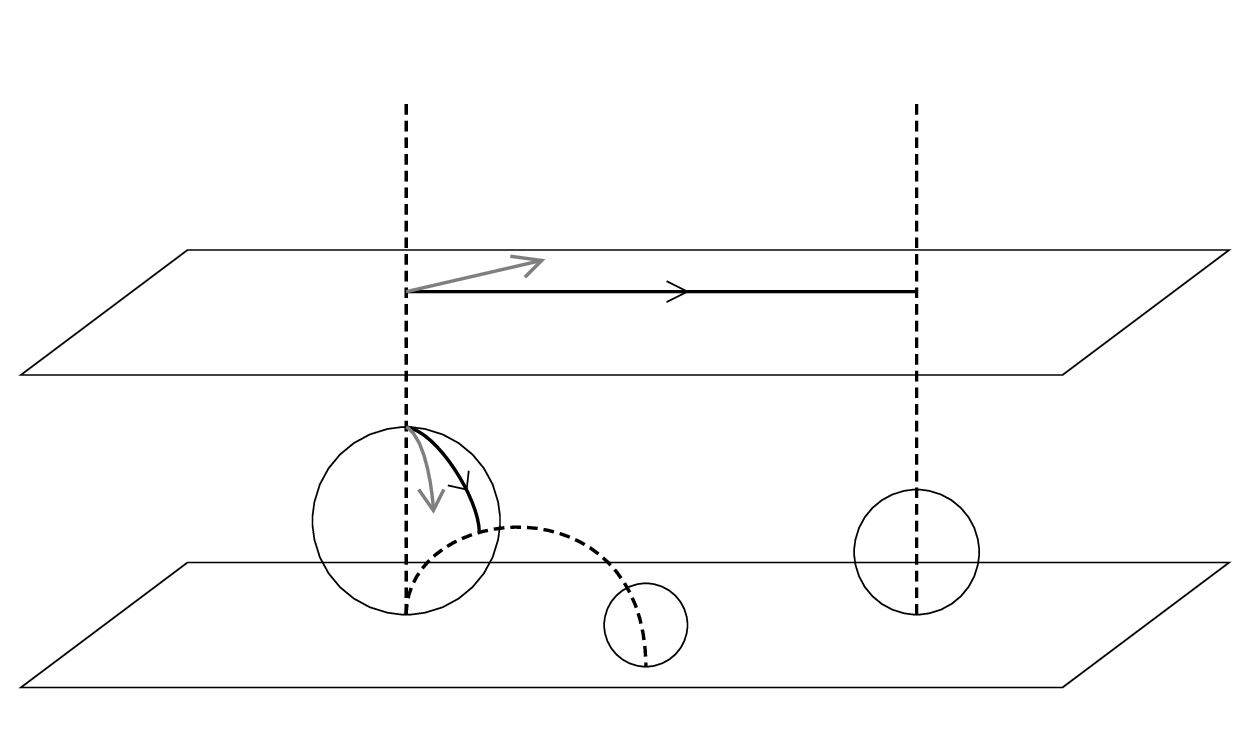}
\caption{}
\end{figure}

It remains to check that the arguments of the complex numbers on each side of equation
{\bf (3)} match.  Fig.\ha 6 illustrates
the situation where the angle \( \phi_{\gamma_1} \) between meridian vectors
\( \mu_1 \bc \mu_2 \) is approximately \( +\pi/2 \) and the arguments of \( u_1 \bc u_2 \)
are each positive acute angles.  The ideal vertices \( z_0 \bc z_1 \bc z_2 \) are situated
at \( |u_1| \bc \infty \bc 0 \) respectively, so \( \arg(\zeta_{\gamma_1}) = \arg(z_3) \).
Noting that the Euclidean ray from \( z_2 \) towards \( z_0 \) is in the positive real
direction, and taking account of the angles between the four vertical planes through \( \gamma_1 \)
containing the vectors \( \mu_1, \mu_2, u_1, u_2 \), we see that
\[ \arg(z_3) \;=\; \arg(u_1) + \arg(u_2) - \phi_{\gamma_1} \;=\; \arg(u_1) + \arg(u_2) - \arg(-w) \h.\]
Therefore, in the notation of the statement of the Proposition, we have
\[ \arg(\overline{\zeta_\gamma}) = -\arg(z_3) = \arg(-w) - (\arg(u_R^{E_1}) + \arg(u_R^{E_2})) \h,\]
and the proof is complete.
\end{proof}

The next proposition gives an alternative version of the equation
\[ \prod_{i=n}^{1}\, \left[\begin{array}{cc}0&-\zeta_i\\1&-1\end{array}\right] \;\sim\; I \]
in M\"{o}bius transformations for a region with \( n \geq 3 \) sides, given by Equation {\bf (1)} of \S 2.
This will be useful for the discussion regarding holonomy representations in \S 5.
\begin{prop}

Let \( R \) be a region of a link diagram with \( n \geq 3 \)
sides, and, starting from some crossing of \( R \), let
\( u_1 \bc w_1 \bc u_2 \bc w_2 \bc \dots \bc u_n \bc w_n \) be the
alternating sequence of edge and crossing labels for \( R \) encountered as one travels
around the region.  Also, for \( 1 \leq i \leq n \) let \( \epsilon_i = 1 \)
(resp. \( \epsilon_i = -1 \)) if the direction of the edge corresponding to \( u_i \)
is with (resp. against) the direction of travel.  Then the equation in M\"{o}bius
transformations for \( R \) may be written as
\begin{equation}
{\text{$\prod_{i=n}^{1}\,$}}
\left(
\left[\begin{array}{cc}0&-w_i\\1&0\end{array}\right]
\left[\begin{array}{cc}1&\epsilon_i u_i\\0&1\end{array}\right]
\right)
\quad\sim\quad
\left[\begin{array}{cc}1&0\\0&1\end{array}\right]
\quad.
\end{equation}
\end{prop}
\begin{proof}[Proof\\]

The conclusion follows from {\bf (1)} in \S 2, together with
\[
\left[\begin{array}{cc}0&-w_i\\1&0\end{array}\right]
\left[\begin{array}{cc}1&\epsilon_i u_i\\0&1\end{array}\right]
\quad=\quad
\left[\begin{array}{cc}0&-w_i\\1&\epsilon_i u_i\end{array}\right]
\quad=\quad
(-\epsilon_i u_i) \h T_{i+1} \, M_i \, T_i^{-1} \h,
\]
where
\[
T_i = \left[\begin{array}{cc}-\epsilon_{i}u_{i}&0\\0&1\end{array}\right] \quad \text{and} \quad
M_i = \left[\begin{array}{cc}0&-\frac{w_i}{(\epsilon_i u_i)(\epsilon_{i+1} u_{i+1})}\\1&-1\end{array}\right]
    = \left[\begin{array}{cc}0&-\overline{\zeta_i}\\1&-1\end{array}\right]
\]
\end{proof}

Suppose now that we have a reduced alternating diagram with \( c \) crossings.
Attached to the diagram are \( 4c \) edge labels (two for each edge), and \( c \)
crossing labels, making altogether \( 5c \) unknowns for our system of equations.
The diagram has \( c + 2 \) regions, giving rise to \( 3(c + 2) \) ``region''
equations in the labels.  Together with the \( 2c \) ``edge'' equations relating labels on
the two sides of an edge, we have in total \( 5c + 6 \) equations in \( 5c \)
unknowns.  In the next section we examine the relationship between solutions to the equations
and representations of the link group into \( \mathrm{PSL_2(\mathbb C)} \).
\end{section}

\begin{section}{The holonomy representation associated to a solution of the label equations}

The connection between solutions to the label equations and parabolic representations of the
fundamental group of the link complement into \( \mathrm{PSL_2(\mathbb C)} \) is not
quite immediate, but can be established without undue difficulty from first principles,
using the classical Wirtinger presentation.

In the traditional picture of the Wirtinger presentation of an oriented link, one takes a
diagram of the link that resides in the projection plane apart from vertical perturbations
within small neighbourhoods of crossings, and one chooses as basepoint of the link
complement some point above the diagram.  The generators of the fundamental group of
the link complement are then path homotopy classes of loops of form
\( \alpha_i * \mu_i * \overline{\alpha_i} \), where the path \( \alpha_i \) travels in a
straight line
from the basepoint down to a point a small distance above the $i$th overpass of the link,
\( \mu_i \) is a small circular loop bounding a disk punctured by the overpass (its direction
being determined by the right-hand screw convention), and \( \overline{\alpha_i} \) is
the reverse of \( \alpha_i \).

For the current context, we shall take the tree that is the union of the \( \alpha_i \),
and push it, keeping the terminal ends of the \( \alpha_i \) fixed, so that the basepoint
lies on the peripheral torus, directly above the
$i$th overpass, and so that each path \( \alpha_i \) is a succession of subpaths, each of
which either travels along the peripheral torus between crossings or travels along the
intercusp segments of crossing geodesics.  At this stage the endpoints of the \( \alpha_i \),
including the basepoint, are all at the tops of meridional circles on the peripheral torus.
It is convenient for purposes of visualization, however, to make a final adjustment to
the tree, whereby we drag the endpoints of the \( \alpha_i \) halfway around their respective
meridional circles so that each lies {\it underneath} its overpass, on the vertical arc that
joins overpass to underpass.

One can imagine a large playground construction consisting of knotted tubes of
cross-sectional diameter say 1 metre, and ladders joining underpasses
with overpasses; the task of \( \alpha_i \) is to travel from one overpass to another
by clambering along sections of tube and climbing or descending ladders, in such a way that
the union of the \( \alpha_i \) is the result of a deformation of the ``classical'' tree
as described above.  We note that if we take the first overpass to be that where the
basepoint is situated, then we may take \( \alpha_1 \) to be the trivial path at
the basepoint.  Also, we may assume that each subpath of \( \alpha_i \) that lies on
the peripheral torus is either of the type \( \alpha_R^E \) of \S 3 (in which case
it corresponds to an edge label), or else is a meridional circle.  The reason why meridional
circles might be needed is that one is not permitted to travel along an underpass under an overpass,
so in order to get past the overpass one would need to climb up to it, travel around a meridional
circle and then climb back down before proceeding.  Naturally, each subpath that
travels along a vertical arc between overpass and underpass corresponds to a crossing label.

\begin{figure}[ht!]
\labellist
\small
\pinlabel $\mathrm{\alpha_{1}}$ at -9 110
\pinlabel $\mathrm{x_{0}}$ at 53 193
\pinlabel $\mathrm{\alpha_{2}}$ at 97 94
\pinlabel $\mathrm{\alpha_{3}}$ at 147 160
\pinlabel $\mathrm{x_{0}}$ at 508 24
\pinlabel $\mathrm{\alpha_{2}}$ at 627 -10
\pinlabel $\mathrm{\alpha_{3}}$ at 639 65
\endlabellist
\centering

\includegraphics[scale=0.36]{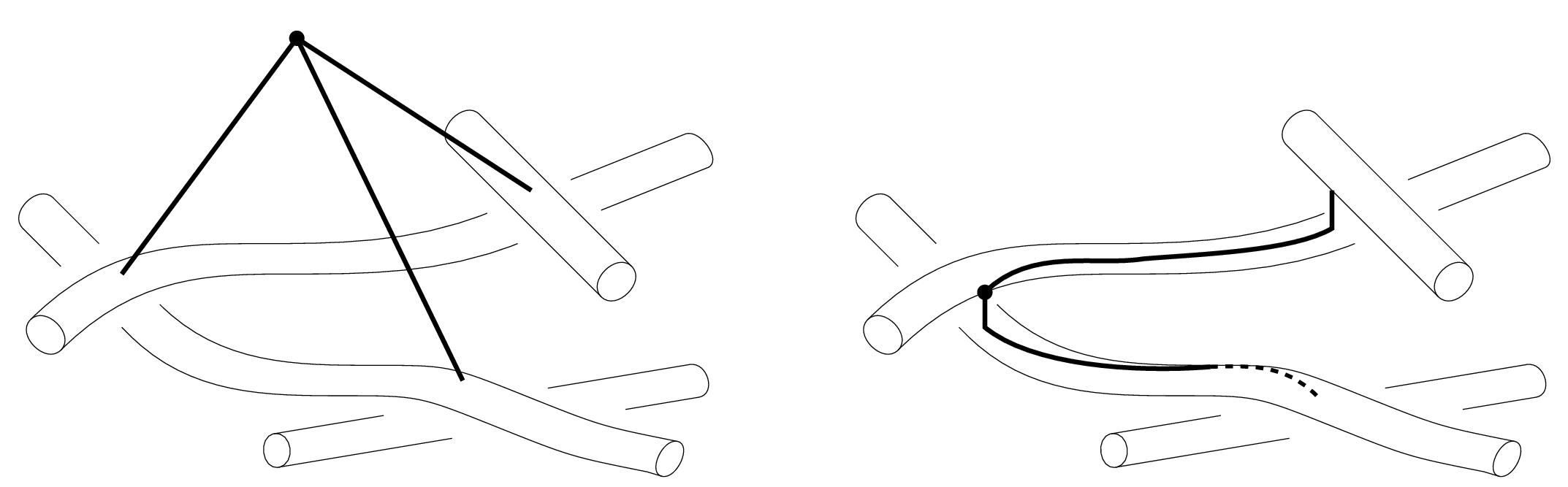}
\caption{}
\end{figure}

The process of deforming the tree \( \bigcup \alpha_i \) is illustrated in Fig.\ha 7.
After the deformation, \( \alpha_1 \) has been shrunk to a point, \( \alpha_2 \) descends to
the underpass at the left-hand crossing before travelling along the peripheral torus
to its target overpass, and \( \alpha_3 \) travels along the peripheral torus from the basepoint
\( x_0 \) before ascending vertically to its target overpass.

We are now ready to specify the representation
\( \phi: \pi_1(S^3 - L) \to \mathrm{PSL_2(\mathbb C)} \) associated to a set of labels
satisfying the label equations.  Each Wirtinger generator
\( [\alpha_i * \mu_i * \overline{\alpha_i}] \) will map to a conjugate of
the parabolic {\( \left[\begin{array}{cc}1&1\\0&1\end{array}\right] \)}, by a matrix \( M_i \)
determined by the path \( \alpha_i \).  Evidently there are choices involved in how
the ``classical'' tree is pushed to its present form; one can travel in either of
two ways around a region between two crossings incident to that region; however, it will be seen
from the definition of \( M_i \) given shortly that independence of those choices is guaranteed
by the label equations, in the form given in {\bf (4)} of Proposition 4.2.

Let us write \( \alpha_i \) as a concatenation of subpaths
\( \alpha_i^1 * \alpha_i^2 * \dots * \alpha_i^{n_i} \), where each subpath is of one of the
three types (i) corresponds to an edge label, (ii) is a meridional circle, (iii) corresponds
to a crossing label.  The standard process of lifting the loop
\( \alpha_i * \mu_i * \overline{\alpha_i} \) to \( \mathbb H^3 \) dictates that we define
\( M_i \;=\; M_i^1 \, M_i^2 \, \dots \, M_i^{n_i} \), where

{
\[ M_{\, i}^{\, j} \;=\; \begin{cases}
\phn\left[\begin{array}{cc}1&\pm u \phl \\0&1\end{array}\right]
      \quad\text{\parbox[t]{3.2in}{\small\rm($\alpha_i^j$ corresponds to the edge label $u$, the
      sign being \( +1 \) iff the direction of the edge agrees with that of \( \alpha_i \))}}\\[3ex]
\phn\left[\begin{array}{cc}1&\pm 1\phl\\0&1\end{array}\right]
      \quad\text{\parbox[t]{3.2in}{\small\rm($\alpha_i^j$ is a meridional circle, the sign
      determined by the RH screw rule)}}\\[3ex]
\phn\left[\begin{array}{cc}0&-w\\1&0\end{array}\right]
      \quad\text{\small($\alpha_i^j$ corresponds to the crossing label $w$)}
\end{cases}
\]
}

Finally, the image of the $i$th Wirtinger generator is given by
\[ \phi([\alpha_i * \mu_i * \overline{\alpha_i}] \h=\h
      M_i \; \left[\begin{array}{cc}1&\phl 1 \phl\\0&1\end{array}\right] \; (M_i)^{-1} \h.\]

\begin{figure}[ht!]
\labellist
\small
\pinlabel $\mathrm{1}$ at 140 143
\pinlabel $\mathrm{2}$ at 377 200
\pinlabel $\mathrm{3}$ at 139 259
\pinlabel $\mathrm{4}$ at 197 22
\pinlabel $\mathrm{5}$ at 252 259
\pinlabel $\mathrm{6}$ at 19 200
\pinlabel $\mathrm{7}$ at 251 143
\pinlabel $\mathrm{8}$ at 197 380
\endlabellist
\centering

\includegraphics[scale=0.45]{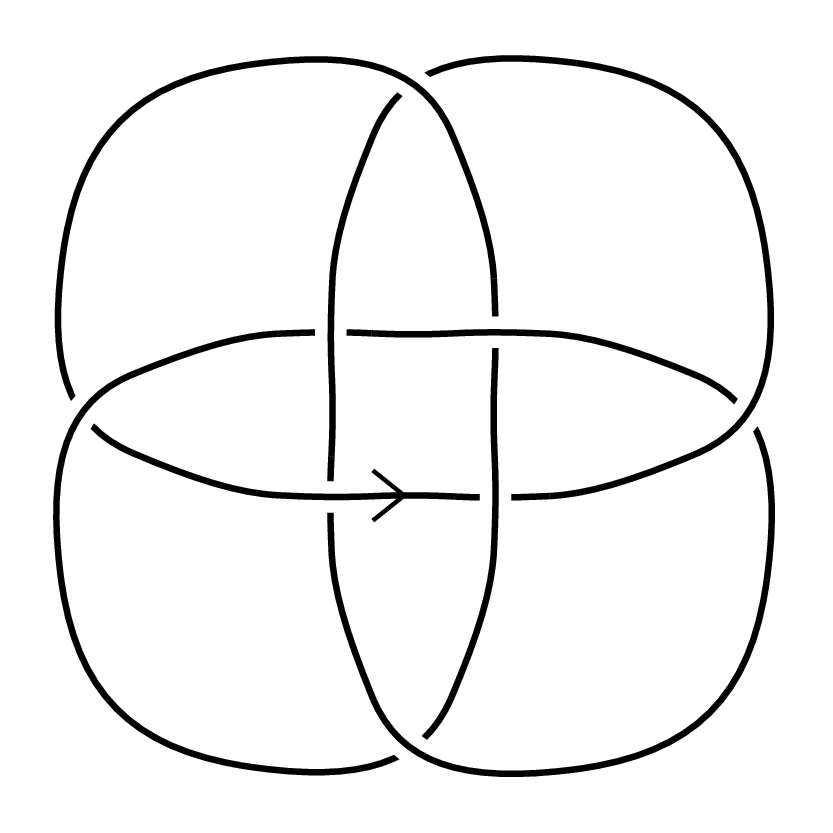}
\caption{}
\end{figure}

There now follows an illustration of the construction of the representation \( \phi \).
The solutions of the label equations for the Turk's head knot (Fig.\ha 8) are obtained
in Example 6.2 below, in the case \( n = 4 \).  As explained in that example, there are two geometric
solutions (forming a complex conjugate pair) and two real non-geometric solutions.

In Fig.\ha 8, the overpasses are labelled from 1 to 8, in order around the knot.  From the knot's
symmetry, the four edge labels for the central 4-sided region are equal, as are that region's
four crossing labels.  Let us denote the common edge label for this region \( u \), and the
common crossing label \( w \).  A quick study of the Wirtinger presentation for this
diagram reveals that the knot group is generated by the Wirtinger generators at
overpasses labelled \( 1, 7, 5 \); let us denote these generators \( a \bc b \bc c \) respectively,
and let us use upper-case \( A \bc B \bc C \) to denote their inverses.
The defining relators for the group (also obtained from the Wirtinger presentation) are then
\begin{align*}
r_1 \;&=\; BabACbcBabCBcaBAbC \h,\\
r_2 \;&=\; ACbcBabCBcaCbcBAbCBcbcBabCBcACbcBAbCBcaB \h;
\end{align*}
we could of course have had much shorter relators at the expense of a larger generating set.

If we put the basepoint at overpass 1, then we can get to overpass 7 by travelling along the
part of the peripheral torus corresponding to the lower edge of the central 4-sided region, 
followed by a climb from the underpass to the overpass at the crossing labelled 7.
In order to get to overpass 5, we have to repeat this process, travelling along the right-hand
edge of the central 4-sided region.
The conjugating matrices for the generators \( b \bc c \) are therefore
\[ M_b \;=\; \left[\begin{array}{cc}1&u\\0&1\end{array}\right] \; \left[\begin{array}{cc}0&-w\\1&0\end{array}\right]
\quad,\quad M_c \;=\; M_b^{\;2} \h.\]
From the computation in Example 6.2, the values of \( u \) for each of the four solutions are:
\[\frac{-1 + i\,\sqrt{4\sqrt{2} - 5}}{\sqrt{2}} \h,\h \frac{-1 - i\,\sqrt{4\sqrt{2} - 5}}{\sqrt{2}} \h,\h
\frac{1 + \sqrt{4\sqrt{2} + 5}}{\sqrt{2}} \h,\h \frac{1 - \sqrt{4\sqrt{2} + 5}}{\sqrt{2}} \h,\]
these being the roots of the polynomial \( x^4 - 6x^2 - 8x - 4 \).
It is evident from the regularity of the central 4-sided region that in each case \( w = \frac{1}{2} u^2 \).
Also, from our choice of basepoint, in each case \( \phi(a) = \left[\begin{array}{cc}1&1\\0&1\end{array}\right]\).
For the first solution, we find that
{\small
\[ \phi(b) \;=\; \left[\begin{array}{cc} 
\frac{1}{\sqrt{2}}\left(2\sqrt{2} + 1 + i\,\sqrt{2\sqrt{2} + 1}\right) & 2\\
\frac{1}{2}\left(-1 - i\,\sqrt{11 + 8\sqrt{2}}\right) & \frac{1}{\sqrt{2}}\left(-1 - i\,\sqrt{2\sqrt{2} + 1} \right)
\end{array}\right] \]
\[ \phi(c) \;=\; \left[\begin{array}{cc}
\frac{1}{\sqrt{2}}\left(2\sqrt{2} + 1 + i\,\sqrt{2\sqrt{2} + 1}\right) & 1\\
-1 - i\,\sqrt{11 + 8\sqrt{2}} & \frac{1}{\sqrt{2}}\left(-1 - i\,\sqrt{2\sqrt{2} + 1} \right)
\end{array}\right] \quad, \]
}

and for the (non-geometric) third solution
{\small
\[ \phi(b) \;=\; \left[\begin{array}{cc} 
\frac{1}{\sqrt{2}}\left(2\sqrt{2} - 1 - \sqrt{2\sqrt{2} - 1}\right) & 2\\
\frac{1}{2}\left(-1 + \sqrt{-11 + 8\sqrt{2}}\right) & \frac{1}{\sqrt{2}}\left(1 + \sqrt{2\sqrt{2} - 1} \right)
\end{array}\right] \]
\[ \phi(c) \;=\; \left[\begin{array}{cc}
\frac{1}{\sqrt{2}}\left(2\sqrt{2} - 1 - \sqrt{2\sqrt{2} - 1}\right) & 1\\
-1 + \sqrt{-11 + 8\sqrt{2}} & \frac{1}{\sqrt{2}}\left(1 + \sqrt{2\sqrt{2} - 1}\right)
\end{array}\right] \quad. \]
}

Representations corresponding to the other two solutions may be constructed similarly.
That all four representations satisfy the relations \( r_1 = 1 \bc r_2 = 1 \) may be verified using software
such as Maple or Mathematica.

Clearly the induced parabolic representation of the link group into \( \mathrm{PSL_2(\mathbb C)} \) is a continuous
function of the edge and crossing labels; also, a variation of the labels within the solution space would
change the geometry of the link complement.  Therefore from Mostow-Prasad rigidity we have:

\begin{thm}
The solution of the label equations corresponding to the geometric structure is isolated.
\end{thm}

\begin{conj}
For hyperbolic alternating links, the solution space of the label equations
is $0$--dimensional, whence the same is true of the space of parabolic representations of the link group into
\( \mathrm{PSL_2(\mathbb C)} \).
\end{conj}

Since the label equations have a reasonably pleasant form for alternating diagrams, it is
reasonable to hope that a proof of Theorem 5.1, and indeed of Conjecture 5.2, can be found without
recourse to Mostow-Prasad rigidity.

We do not know of a simple test for deciding which solution of the equations
corresponds to the geometric structure; however, by subdividing the link complement into ideal
tetrahedra and keeping track of labels, one may compute
the volume of the representation given by the solution, and the solution with the greatest
volume will be the geometric solution \cite{Franc}.  Empirically, for alternating links the
geometric solution is that for which the ``region'' ideal polygons are closest to
being regular.

Experimentally, the near-regularity of these polygons is particularly
evident for alternating diagrams of so-called Conway basic polyhedra, these being links
possessing alternating diagrams with no two-sided region.  The simplest examples are
(i) the Borromean rings, (ii) the Turk's head knot, (iii) the 9-crossing knot denoted 9* by
J.H. Conway in \cite{Con} (listed as $9_{40}$ in \cite{Rolf} and as $9a37$ in
the Dowker-Thistlethwaite classification).

For the Borromean rings, the regions are forced to
be regular by virtue of being 3-sided, and regularity of the 4-sided regions
of the Turk's head knot is a consequence of the symmetries of that knot; see Example 6.2
below.  On the other hand, as already seen in Example 2.3.1, the knot 9* (Fig. 1)
has 4-sided regions that deviate slightly from being regular.
\end{section}

\begin{section}{Examples}

In this section we give some examples illustrating the use of the hyperbolicity
equations from \S 4.

\begin{subsection}{\rm The figure-eight knot}

\begin{figure}[ht!]
\labellist
\small
\pinlabel $\mathrm{u_{1}}$ at -169 -136
\pinlabel $\mathrm{u_{2}}$ at -50 -68
\pinlabel $\mathrm{u_{3}}$ at 49 -70
\pinlabel $\mathrm{u_{4}}$ at 169 -136
\pinlabel $\mathrm{u_{1} + 1}$ at -120 -135
\pinlabel $\mathrm{u_{2} + 1}$ at -93 -86
\pinlabel $\mathrm{u_{3} + 1}$ at 93 -86
\pinlabel $\mathrm{u_{4} + 1}$ at 120 -135
\pinlabel $\mathrm{0}$ at 1 10
\pinlabel $\mathrm{0}$ at 1 -11
\pinlabel $\mathrm{0}$ at -12 -177
\pinlabel $\mathrm{0}$ at 12 -177
\pinlabel $\mathrm{1}$ at -38 -177
\pinlabel $\mathrm{1}$ at 38 -177
\pinlabel $\mathrm{-1}$ at -1 39
\pinlabel $\mathrm{-1}$ at -1 -39
\pinlabel $\mathrm{w_{1}}$ at -71 1
\pinlabel $\mathrm{w_{1}}$ at 71 1
\pinlabel $\mathrm{w_{2}}$ at 0 -115
\pinlabel $\mathrm{w_{2}}$ at 0 -246
\pinlabel {\large $\mathrm{I}$} at -86 -194
\pinlabel {\large $\mathrm{II}$} at 0 -77
\pinlabel {\large $\mathrm{III}$} at 86 -194
\endlabellist
\centering

\includegraphics[scale=0.65]{fig8}
\caption{}
\end{figure}

Recall that for a 3-sided region all shape parameters are 1.  Regions I, II, III each
provide three equations in the labels, as follows:
\begin{align*}
\mathrm{I{:}\phantom{II}}\quad &w_2 = -(u_2+1) \h,\h w_2 = -(u_1+1) \h,\h w_1 = (u_1+1)(u_2+1)\\
\mathrm{II{:}\phantom{I}}\quad &w_{1}=u_{2} \h,\h w_1=u_{3} \h,\h w_{2}=u_{2}u_{3}\\
\mathrm{III{:}          }\quad &w_2 = -(u_3+1) \h,\h w_2 = -(u_4+1) \h,\h w_1 = (u_3+1)(u_4+1)
\end{align*}

Collecting these results, we obtain
\[ w_{1}=u_{1}=u_{2}=u_{3}=u_{4} \h,\h u_{1}^2+u_{1}+1=0 \h,\h
w_{2}=-(u_{1}+1) \h,\h w_{1}=w_{2}^{2} \h.\]
Therefore, without loss of generality, we have
\begin{align*}
u_{1}=u_{2}=u_{3}=u_{4} &= \frac{1}{2}(-1 + i\sqrt{3})\h,\\
w_{1} &= \frac{1}{2}(-1 + i\sqrt{3})\h,\\
w_{2} &= \frac{1}{2}(-1 - i\sqrt{3})\h.
\end{align*}
\end{subsection}

\begin{subsection}{\rm The closure \( L_n \) of the braid
\( (\sigma_1 \sigma_2^{-1})^n \), \( n \geq 3 \)}

\begin{figure}[ht!]
\labellist
\pinlabel $\mathrm{u_{1}}$ at 122 142
\pinlabel $\mathrm{u_{1}}$ at 122 -56
\pinlabel $\mathrm{u_{2}}$ at -22 110
\pinlabel $\mathrm{u_{2}}$ at -22 9
\pinlabel $\mathrm{\overline{u_{1}}}$ at -124 58
\pinlabel $\mathrm{w}$ at 50 77
\pinlabel $\mathrm{w}$ at 50 -122
\pinlabel $\mathrm{\overline{w}}$ at -50 -75
\pinlabel {\normalsize (i) A section of the braid} at 0 -250
\pinlabel {\normalsize (ii) The closure of $\left( \sigma_1 \sigma_2^{-1} \right)^8$} at 500 -250
\endlabellist
\centering

\includegraphics[scale=0.43]{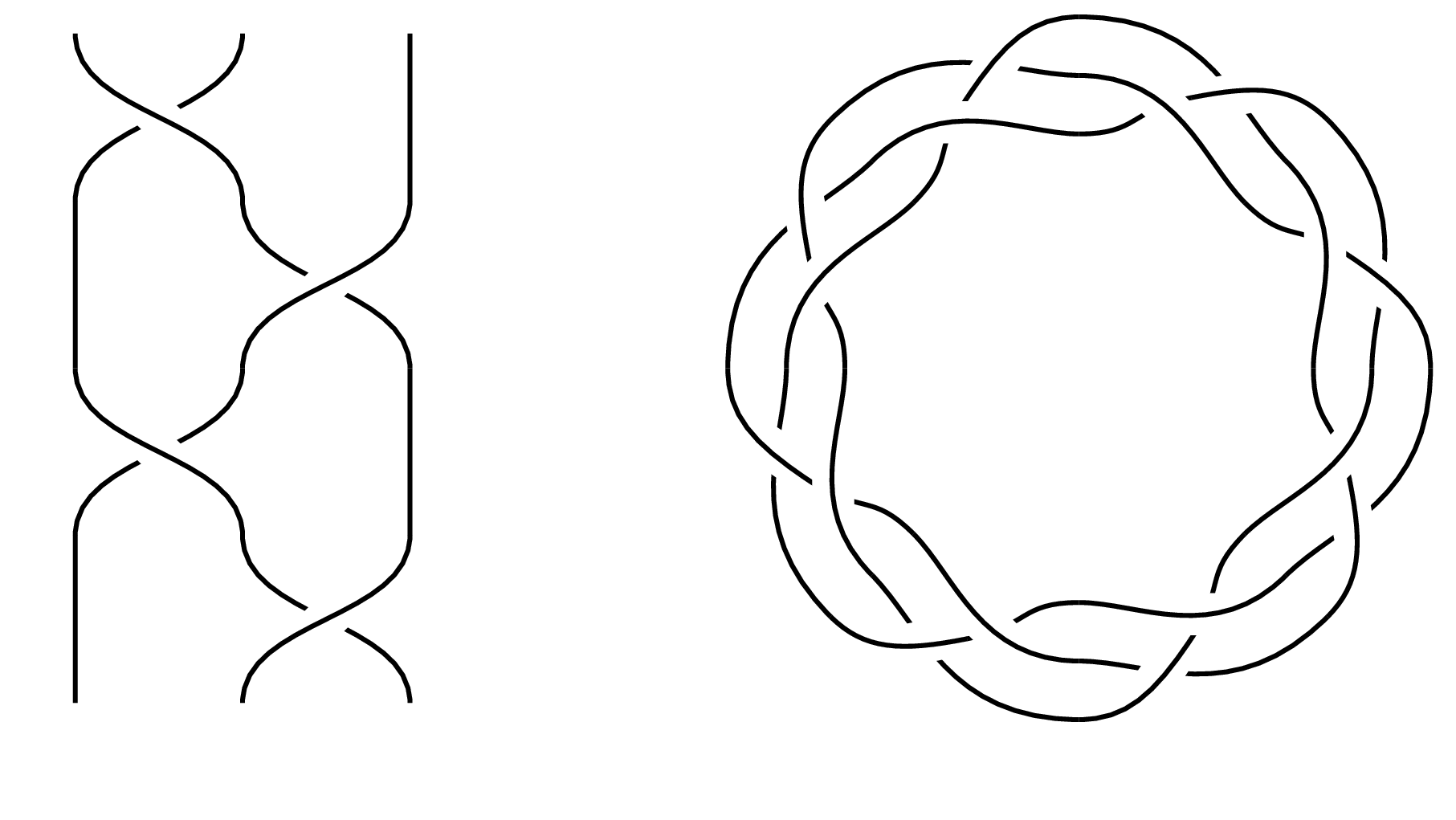}
\caption{}
\end{figure}

In Fig.\ha 10(i) we have exploited symmetries of the link \( L_n \) in order to economize on labels.
For notational convenience let \( \lambda_n \) denote \(\frac{1}{2} \sec \frac{\pi}{n} \),
the positive square root of the shape parameter of a regular ideal $n$-gon.
Looking at Fig.\ha 10(i), from the right-hand $n$--sided region we have
\[ w = \lambda_n^2\, u_1^2 \h,\]
and from 3--sided regions we have
\[ w = -(u_1+1)(u_2+1) \h,\h w = u_2^2 \h.\]
Thus \( u_1 = \pm\frac{1}{\lambda_n} u_2 \); however, for the geometric solution we must take
\( u_1 = \frac{1}{\lambda_n} u_2 \), as otherwise all edge labels will turn out to be real, resulting in
a collapse of the peripheral structure.  It then follows quickly that \( u_2 \) 
satisfies
\[ \left(1 + 2\cos\frac{\pi}{n}\right)u_2^2 + \left(1 + 2\cos\frac{\pi}{n}\right)u_2 + 1 = 0 \h,\]
the two solutions of this quadratic yielding geometric structures corresponding to
the two orientations of the link complement.

If \( n \) is divisible by \( 3 \), \( L_n \) is a link of \( 3 \) components, with symmetries
acting transitively on the set of components; otherwise \( L_n \) is a knot.
The link \( L_3 \) is the Borromean rings, and \( L_4 \) is the Turk's head knot.  In the case
of the Borromean rings, the crossing labels are \( \pm\frac{i}{2} \), indicating that the
overpass and underpass at each crossing are perpendicular to one another; this is
also evident from the extra symmetries possessed by this link.

It is of some interest to note that as \( n \to \infty \), \( |w| \to \frac{1}{3} \).
This has the following interpretation in terms of meridian lengths.
Let us expand cusp cross-sections, keeping meridian lengths equal, until their
union just ceases to be embedded in the link complement (in the case where \( L_n \) is
a knot the cusp will touch itself).  Let \( \ell_n \) denote the length of a meridian on
one of these expanded cusp cross-sections; then \( \ell_n \) tends to a limit \( \sqrt{3} \)
as \( n \to \infty \).

Another interesting aspect of the links \( L_n \) concerns the canonical cell decompositions\\
\cite{EP}
of their complements.  The alternating diagram of \( L_n \) exhibits a decomposition of
the link complement into two congruent ideal polyhedra, one ``above'' and one ``below'' \cite{AR},
and this decomposition is precisely the canonical cell decomposition, as is readily verified
by examination of the horoball pattern.  See also the discussion of canonical cell decompositions
in \cite{SW}.

The ideal polygons corresponding to the regions of the alternating diagram of \( L_n \) are
all regular, hence also planar (compare with the knot 9*, discussed in Example 2.3.1)\,.
\end{subsection}

\begin{subsection}{\rm Three-punctured sphere}

Let us consider the three-punctured sphere \( S \) in the configuration of Fig.\ha 11, with
meridional punctures at the two parallel vertical strands, and with a longitudinal puncture at
the circular link component.  It is well-known \cite{Ad3} that we may take \( S \)
to be a totally geodesic surface, constructed by gluing two ideal triangles together along
their edges.  If we choose cusps of an ideal triangle so that the length of each cusp
boundary is \( 1 \), then these cusp boundaries will be tangent to one another; therefore
if we choose cusps of a three-punctured sphere so that each cusp boundary is a circle
of length \( 1+1 = 2 \), these circles will touch one another.  If we now retract the
cusps of the three-punctured sphere so that their boundaries have length \( 1 \),
the intercusp length of each of the three geodesics joining distinct punctures will be
\( \log(4) \).  Since by convention our meridians always have length \( 1 \), it follows
that \( w_{3} \), the label for the geodesic represented by the horizontal line at the top of Fig. 11,
has modulus equal to\( \frac{1}{4} \) (recall from the remark following Conjecture 3.2 that
to any geodesic joining cusps there is an associated complex number, defined exactly
as crossing labels are defined).  On the other hand, since the three-punctured
sphere is totally geodesic, the meridians belonging to the vertical strands lie in the
same hyperbolic plane as this geodesic, whence \( w_{3} = \frac{1}{4} \).
({\it Note.}  Here we have chosen orientations of the strands so that they are
parallel; were they anti-parallel, instead we would have \( w_{3} = -\frac{1}{4} \).)

\begin{figure}[ht!]
\labellist
\pinlabel $\mathrm{w_{1}}$ at -102 95
\pinlabel $\mathrm{w_{1}}$ at -102 -98
\pinlabel $\mathrm{w_{2}}$ at 102 95
\pinlabel $\mathrm{w_{2}}$ at 102 -98
\pinlabel $\mathrm{w_{3}}$ at 0 180
\pinlabel $\mathrm{u_{1}}$ at 0 55
\pinlabel $\mathrm{u_{2}}$ at 0 -55
\pinlabel $\mathrm{1}$ at -55 0
\pinlabel $\mathrm{-1}$ at 55 0
\pinlabel {\Large $\mathrm{I}$} at 0 0
\pinlabel {\Large $\mathrm{II}$} at 0 120
\endlabellist
\centering

\includegraphics[scale=0.50]{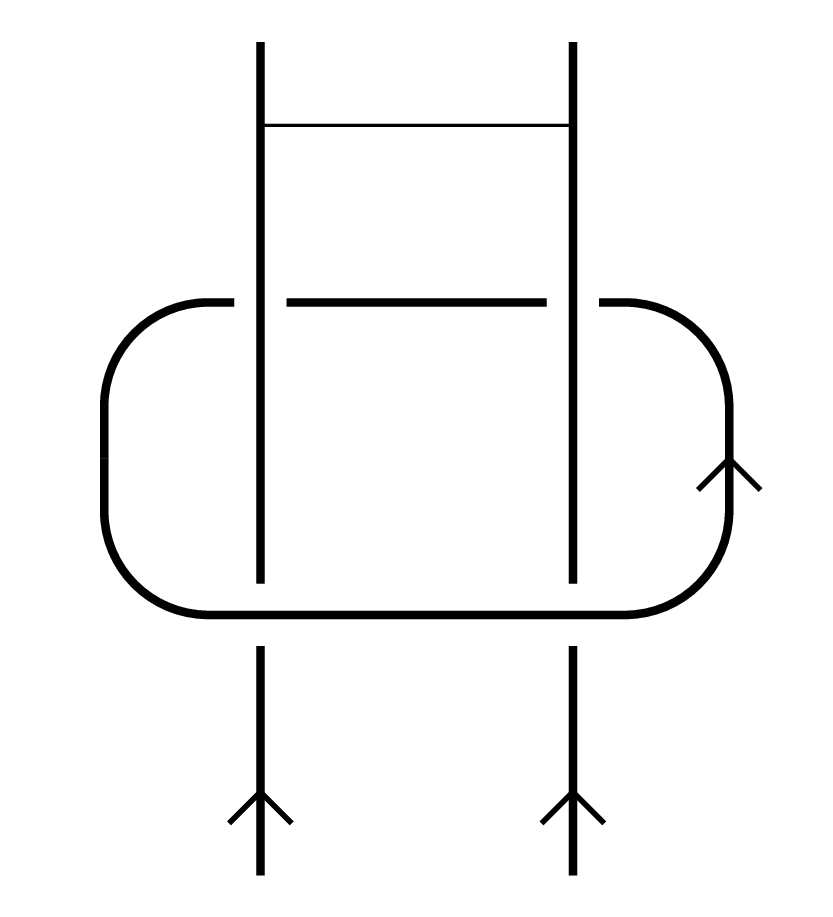}
\caption{}
\end{figure}

From the 4-sided region marked I, from the relation \( f_4 \) in shape parameters given
in \S 2, we have
\[ -\frac{w_{1}}{u_{2}}-\frac{w_{1}}{u_{1}}=1 \h,\h
   -\frac{w_{1}}{u_{1}}+\frac{w_{2}}{u_{1}}=1 \h,\h
    \frac{w_{2}}{u_{1}}+\frac{w_{2}}{u_{2}}=1 \h,\]
from which it follows easily that
\[ w_{2} = -w_{1} \quad \text{and} \quad u_{1} = u_{2} = -2w_{1} \h.\]
Then, from the 3-sided region marked II, we have
\[ w_{3} = \frac{-w_{1}w_{2}}{u_{1}^{2}} = \frac{1}{4} \h,\]
thus recovering the result obtained by the geometric argument of the previous paragraph.
\end{subsection}
\end{section}

\begin{section}{Labels on tangles}

Let \( L \) be a hyperbolic link in \( S^3 \).  In this section we are concerned with
tangles \( (B \bc T) \subset (S^3 \bc L) \) , where \( B \) is a $3$--ball and \( T \)
is a properly embedded $1$--dimensional submanifold of \( B \) meeting the boundary
of \( B \) transversely in four points.  Thus \( \partial B - \partial T \) is a
$4$--punctured sphere; we require that this ``Conway sphere'' be essential in \( S^3 - L \).  Thus
both \( (B \bc T) \) and the complementary tangle \( (\overline{S^3 - B} \bc \overline{L - T}) \)
are {\it non-trivial} in the sense of \cite{Lickorish}, {\it i.e.} neither is homemomorphic
as a pair to \( (B_0 \bc T_0) \), where
\( B_0 = \{ (x_1,x_2,x_3) \in \mathbb R^3: x_1^2 + x_2^2 + x_3^2 \leq 1 \} \) and
\( T_0 = \{ (x_1,x_2,x_3) \in B_0: x_2 = 0 \text{\;and\;} x_3 = \pm \frac{1}{2} \} \).
The restriction of the hyperbolic metric on \( S^3 - L \) to \( B - T \) is complete
(in the sense that Cauchy sequences converge), as is the restriction to the $4$--punctured sphere
\( \partial B - \partial T \).

Our first observation is that the boundary of \( (B \bc T) \) enjoys a certain kind of
symmetry, which can be described as follows.  Let us suppose that \( B \) is a standard 3--ball
meeting the projection plane in an equatorial disk \( \Delta \), and that the tangle \( T \) is
contained in this disk except for small vertical perturbations at crossings; thus the four
boundary points of \( T \) are contained in the circle \( C = \partial\Delta \).  Let these points,
taken in cyclic order around \( C \), be \( Q_{1} \bc Q_{2} \bc Q_{3} \bc Q_{4} \).
If our tangle diagram is part of a taut link diagram, then for each \( i \in \{ 1 \bc 2 \bc 3 \bc 4 \} \)
we have a complex number \( w_{i,i+1} \) associated to the sub-arc of \( C \) joining
\( Q_i \) to \( Q_{i+1} \) (suffixes taken modulo 4).

\begin{prop}
In the notation of the previous paragraph, \( w_{1,2} = w_{3,4} \)
and \( w_{2,3} = w_{4,1} \).
\end{prop}

\begin{proof}
This follows essentially from the results of D. Ruberman in \cite{Ruberman}.
In that paper it is shown that if we take the Conway sphere \( \partial B - \partial T \) to
be of least area in its homotopy class, then an involution of mutation is a local isometry near
\( \partial B - \partial T \).  We apply this result to the mutation \( \tau \) corresponding
to a half-turn about an axis perpendicular to the projection plane.  In the notation
of the previous paragraph, let \( \alpha_{i,i+1} \) be the interior of the sub-arc of
\( C \) joining \( Q_i \) to \( Q_{i+1} \); then we may assume that the \( \alpha_{i,i+1} \)
lie on the least area Conway sphere and that the involution \( \tau \) interchanges
neighbourhoods in \( S^3 - L \) of opposite pairs of these arcs.  \( \tau \) also interchanges
intersections of these neighbourhoods with our chosen horospherical cusp boundaries.

It is sufficient to check that \( w_{1,2} = w_{3,4} \), as the other equality will follow
simply by shifting indices.
In order to check that \( w_{1,2} = w_{3,4} \) holds, we need to compare the relative positions
of horospheres at each end of a lift \( \widetilde{\alpha_{1,2}} \) of \( \alpha_{1,2} \),
with those at each end of a lift \( \widetilde{\alpha_{3,4}} \) of \( \alpha_{3,4} \).
Since the involution interchanges open patches of cusp boundaries, it follows that
there is an isometry of \( \mathbb H^3 \) mapping the pair of horospheres at the ends
of \( \widetilde{\alpha_{1,2}} \) to the pair at the ends of \( \widetilde{\alpha_{3,4}} \);
we conclude that \( |w_{1,2}| = |w_{3,4}| \).  Equality of the respective arguments of
\( w_{1,2} \bc w_{3,4} \) follows from consideration of the action of \( \tau \) on the
affine structures on the open patches of cusp boundaries.  Indeed, \( \tau \) maps meridian
curves at \( Q_1 \bc Q_2 \) to meridian curves at \( Q_3 \bc Q_4 \), and since the tangle
has two inward-pointing ends and two outward-pointing ends, \( \tau \) either preserves or
reverses both of their orientations.
\end{proof}

\begin{rem}
The conclusion of Proposition 7.1 holds also for diagrams of
tangles that are trivial in the sense of \cite{Lickorish}.  This can be established
directly from the equations governing edge and crossing labels, by taking the diagram
to be of standard 4-plait ``rational tangle'' form \cite{Con} and proceeding by
induction on the number of crossings of the 4-plait, the inductive step corresponding
to the addition of a crossing by performing a half-twist of two adjacent ends of
the tangle.  The basis for the induction is the conclusion for a tangle diagram
of a single crossing, with crossing label \( w \);
there each arc \( \alpha_{i,i+1} \) is homotopic to the arc travelling vertically from
underpass to overpass at the crossing, whence each \( w_{i,i+1} \) is equal to \( w \).
\end{rem}

\begin{cor}
If a crossing \( c_{1} \) is traded for a crossing \( c_{2} \) by
a flype of the link diagram, then the crossings \( c_{1} \bc c_{2} \) have equal crossing labels.
\end{cor}

\begin{proof}
With appropriate numbering, the crossing labels for \( c_{1} \bc c_{2} \) are
precisely \( w_{1,2} \bc w_{3,4} \) for the tangle turned upside-down by the flype
(see Fig.\ha 12).
\end{proof}

\begin{figure}[ht!]
\labellist
\pinlabel $\mathrm{c_{1}}$ at -54 -28
\pinlabel $\mathrm{c_{2}}$ at 454 -28
\endlabellist
\centering

\includegraphics[scale=0.6]{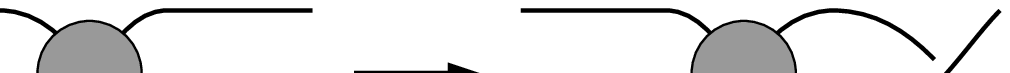}
\caption{A flype:  the tangle represented by the shaded disk is rotated through
a half-turn about a horizontal axis, untwisting the left-hand strands and twisting
the right-hand strands, thus trading the crossing \( c_1 \) for the crossing \( c_2 \).}
\end{figure}

We turn now to consideration of complete hyperbolic structures on a tangle
complement \( B - T \).  Given a solution to the equations in edge and crossing
labels corresponding to the regions of a diagram of the tangle, the method of
\S 5 will yield a corresponding parabolic representation of the fundamental group of
\( B - T \) into \( \mathrm{PSL_2(\mathbb C)} \).  Here we are considering a tangle
pair \( (B \bc T) \) in isolation, not necessarily contained in a link pair
\( (S^3 \bc L) \).

To avoid dealing with manifolds with non-compact boundary (in this case a
4--punctured sphere), the appropriate setting for investigating the Teichm\"{u}ller
space of complete hyperbolic structures on a tangle complement is that of a {\it pared manifold}
\cite{Morgan} (see also \cite{CM}):  one constructs a manifold \( M \), homotopy
equivalent to \( B - T \) and with compact boundary, by excising from \( B \) the
interior of a regular neighbourhood of \( T \); the boundary \( P \subset \partial M \)
of this regular neighbourhood is marked as being peripheral, and is termed the
{\it parabolic locus} of the pared manifold \( (M \bc P) \).  The reader is referred
to Definition 4.8 of \cite{Morgan} for details.

From a discrete faithful parabolic representation of \( \pi_1(B - T) \) into
\( \mathrm{PSL_2(\mathbb C)} \), or equivalently a discrete faithful representation
of \( \pi_1(M) \) into \( \mathrm{PSL_2(\mathbb C)} \) mapping elements of
\( \pi_1(P) \) to parabolics, one obtains a complete hyperbolic structure on the
pared manifold \( (M \bc P) \).  

\begin{prop}
Let \( (M \bc P) \) be a pared manifold arising from
a tangle pair \( (B \bc T) \) as above.  Then (i) the space \( \mathcal T \) of
complete hyperbolic structures on \( (M \bc P) \), if non-empty, is homeomorphic to
\( \mathbb R^2 \); (ii) \( \mathcal T \), regarded as a set of solutions of the label
equations for a given taut diagram of \( (B \bc T) \), is an open subset of the space
of all solutions of the label equations.
\end{prop}

\begin{proof}
It follows from a classical result of L. Bers \cite{Bers} (see
also \cite{Morgan}, Theorem 9.2) that there is a homeomorphism from the space of
complete hyperbolic structures on \( (M \bc P) \) to the Teichm\"{u}ller space of complete,
finite area hyperbolic structures on \( \partial M - P \).  Since \( \partial M - P \)
is homeomorphic to a 4--punctured sphere, part (i) follows from the well-known fact that
the Teichm\"{u}ller space of finite area complete hyperbolic metrics on a 2-sphere with $n$
punctures (\( n \geq 3 \)) is homeomorphic to \( \mathbb R^{2n-6} \). 
Part (ii) follows directly from the Deformation Theorem in \cite{Marden}, \S 6.4; indeed,
Marden's theorem also tells us that \( \mathcal T \) is a
connected, complex analytic manifold of dimension 1, consistent with part (i)
of Proposition 7.2.
\end{proof}

The Kleinian group \( \Gamma \) that is the image of \( \pi_1(M) \) in
\( \mathrm{PSL_2(\mathbb C)} \) can act either on the universal
cover \( \widetilde{M} \), or on the whole of \( \mathbb H^3 \).  In the former case
we obtain a hyperbolic metric on \( M \), and the latter provides a complete hyperbolic metric on
\( \mathbb H^3/\,\Gamma \;\widetilde{=}\; M - \partial M \); these two scenarios can be considered
to be in some sense equivalent.

The following two-part remark is tangential to the main discussion, but may be of
some interest.

\begin{rem}
{\bf (i)}\h\h From Proposition 7.1, it follows that if \( (B \bc T) \subset (S^3 \bc L) \)
for some hyperbolic link \( L \), and if the boundary 4--punctured sphere of \( (B \bc T) \)
is essential in \( S^3 - L \), then the numbers \( w_{i,j} \) for \( T \) arising
from the geometric structure on \( S^3 - L \) obey the symmetry property
\( w_{1,2} = w_{3,4} \;,\; w_{2,3} = w_{4,1} \).  Thus this property of the \( w_{i,j} \)
holds for those points of the Teichm\"{u}ller space \( \mathcal T \) that correspond to
the situation where we sum \( (B \bc T) \) with a non-trivial tangle to form \( (S^3 \bc L) \)
for some hyperbolic link \( L \).  In fact, if 
\( \mathcal T \neq \emptyset \) there are infinitely many such points of \( \mathcal T \),
as one can form infinitely many distinct hyperbolic links \( L \) by summing \( T \) with
an infinite sequence of pairwise distinct non-trivial tangles.  From Marden's Deformation Theorem
\( \mathcal T \) has complex dimension \( 1 \), so \( w_{1,2} \bc w_{3,4} \) cannot be algebraically
independent over \( \mathcal T \).  The fact that \( w_{1,2} = w_{3,4} \) holds for infinitely many
points of \( \mathcal T \) then implies that the equality holds for all points of \( \mathcal T \);
likewise, \( w_{2,3} = w_{4,1} \) throughout \( \mathcal T \).
\vv

{\bf (ii)}\h\h Suppose that we have diagrams of tangles \( T_1 \bc T_2 \), with
respective solutions \( S_1 \bc S_2 \) to the label equations for these diagrams.
The numbers \( w_{i,\,j} \) of Proposition 7.1 are then determined for each diagram.
Suppose that we form a link diagram by summing together \( T_1 \) and \( T_2 \);
Then the solutions \( S_1 \bc S_2 \) can be amalgamated to create a solution to the label
equations for this link if and only if corresponding \( w_{i,\,j} \) for the two tangles agree.
In the case of geometric solutions, we can think of \( w_{1,2} \bc w_{2,3} \) for
each \( T_i \) as defining a complex curve in \( \mathbb C^2 \); the two curves
then meet generically at a point, this point corresponding to the complete structure on the
link complement.
\end{rem}

We turn our attention to a particular kind of tangle, endowed with a surprising rigidity property.

\begin{defn}
A tangle \( T \) conforming to Fig.\ha 13(i) will be called an
{\it encircled tangle}, and a diagram of that type will be called a {\it standard diagram}
of an encircled tangle.  In that figure, the shaded disk represents an arbitrary tangle
\( U \), and \( T \) is the union of \( U \) and a simple closed curve \( C \) that weaves
around the ends of \( U \) in alternating fashion.
\end{defn}

\begin{defn}
Let \( T \) be an encircled tangle as above, with encircling
simple closed curve \( C \).  We assume that \( T \) is represented by a standard diagram,
{\it i.e.} one conforming to Fig.\ha 13(i).  The four crossing labels attached to crossings of
\( C \), and four edge labels attached to the {\it insides} of the edges of \( C \)
will be called {\it boundary labels}, and the remaining labels of \( T \) interior to \( C \)
will be called {\it interior labels.}
\end{defn}

\begin{figure}[ht!]
\labellist
\pinlabel $\mathrm{C}$ at -86 -95
\pinlabel $\mathrm{U}$ at 45 -37
\pinlabel $\mathrm{u_{1}}$ at 412 65
\pinlabel $\mathrm{u_{2}}$ at 292 65
\pinlabel $\mathrm{u_{3}}$ at 292 -65
\pinlabel $\mathrm{u_{4}}$ at 412 -65
\pinlabel $\mathrm{w_{1}}$ at 479 23
\pinlabel $\mathrm{w_{2}}$ at 374 127
\pinlabel $\mathrm{w_{3}}$ at 220 23
\pinlabel $\mathrm{w_{4}}$ at 374 -127
\pinlabel {\small (i)\, An encircled tangle} at 0 -180
\pinlabel {\small (ii)\, Boundary labels} at 350 -180
\endlabellist
\centering

\includegraphics[scale=0.55]{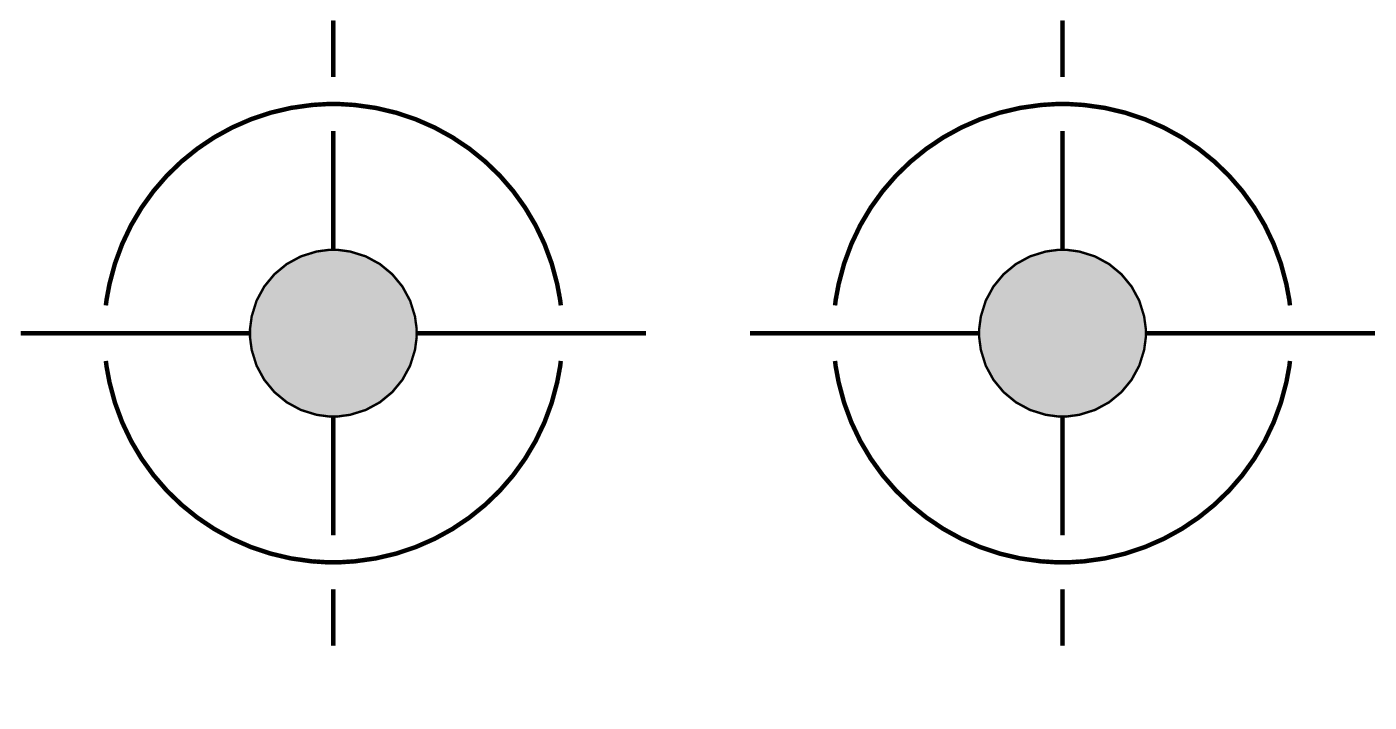}
\caption{}
\end{figure}

The next theorem expresses the rigidity property for encircled tangles.  As before,
we denote by \( \mathcal T \) the space of complete hyperbolic structures on
\( B - T \).  From Proposition 7.2, \( \mathcal T \;\widetilde{=}\; \mathbb R^2 \),
and we may regard \( \mathcal T \) as an open subspace of the solution space of the
label equations.

\begin{thm}

{\bf (i)}\h\h  Each interior label of \( T \) is constant over \( \mathcal T \).
\v

{\bf (ii)}\h\h  Let \( S_1 \bc S_2 \) be solutions of the label
equations corresponding to points of \( \mathcal T \).
There exists a non-zero complex number \( k \) such that if \( \zeta \) is a
boundary label of \( T \) in the solution \( S_1 \), then the corresponding boundary label
in the solution \( S_2 \) is \( k \zeta \).
\end{thm}

\begin{proof}
We assume that the diagram of the subtangle \( U \) of \( T \) (fig. 13)
has at least one crossing; thus each region of the tangle diagram incident to the
encircling curve \( C \) has at least three sides, and shape parameters for these regions
are defined.  The simple case where \( U \) has no crossings will be dealt with directly from
the label equations, in the proof of Theorem 7.4 below.

Recall that the shape parameter at a corner of a region is, up to sign,
the quotient of the label at that crossing by the product of the two incident edge labels.
If we take any non-zero complex number \( k \) and replace each boundary label of \( T \)
in \( S_1 \) by its product with \( k \), by inspection of Fig. 13(ii) all shape
parameters for \( T \) are unchanged, and the equations given by the regions of the tangle
\( T \) are still satisfied.

Let \( \mathcal C \) be the component of the solution set of the label equations that
contains \( \mathcal T \); \( \mathcal C \) has the structure of a complex algebraic
set.  From  A. Marden's  Deformation  Theorem \cite{Marden}, \( \mathcal T \) is open in
\( \mathcal C \), consists only of smooth points of \( \mathcal C \), and has real dimension 2;
it follows that \( \mathcal C \) has complex dimension \( 1 \), {\it i.e.} is a
complex algebraic curve, and that \( \mathcal T \) is away from any singular points
of \( \mathcal C \).  Let \( \mathbf{z}_0 \) be any point of \( \mathcal T \)
(regarded as a subset of \( \mathcal C \)), and let \( \mathcal C' \subset \mathcal C \)
be the set of solutions obtained by taking a non-zero complex number \( k \) and then
multiplying all boundary labels for \( \mathbf{z}_0 \) by \( k \), but keeping interior labels fixed, as in
the previous paragraph (thus \( k \) varies over \( \mathcal C' \)). 
Since \( \mathcal T \subset \mathcal C \) and \( \mathcal C \)
has complex dimension 1,  \( k \) must parametrize the set \( \mathcal T \); in particular,
\( \mathcal T \subset \mathcal C' \).  Both parts of the theorem follow.
\end{proof}

To summarize informally, \( \mathcal T \) has one complex ``degree of freedom'', and because of
dimensional constraints, this ``structural flexibility'' has to be parametrized by \( k \).
Thus the complete structures of an encircled tangle are related in a simple way;
however, we do not know how to determine which values of \( k \), applied in this
manner to the particular complete structure \( \mathbf{z}_0 \), yield points of
\( \mathcal T \).

\begin{figure}[ht!]
\labellist
\pinlabel $\mathrm{\gamma}$ at 657 -30
\endlabellist
\centering

\includegraphics[scale=0.36]{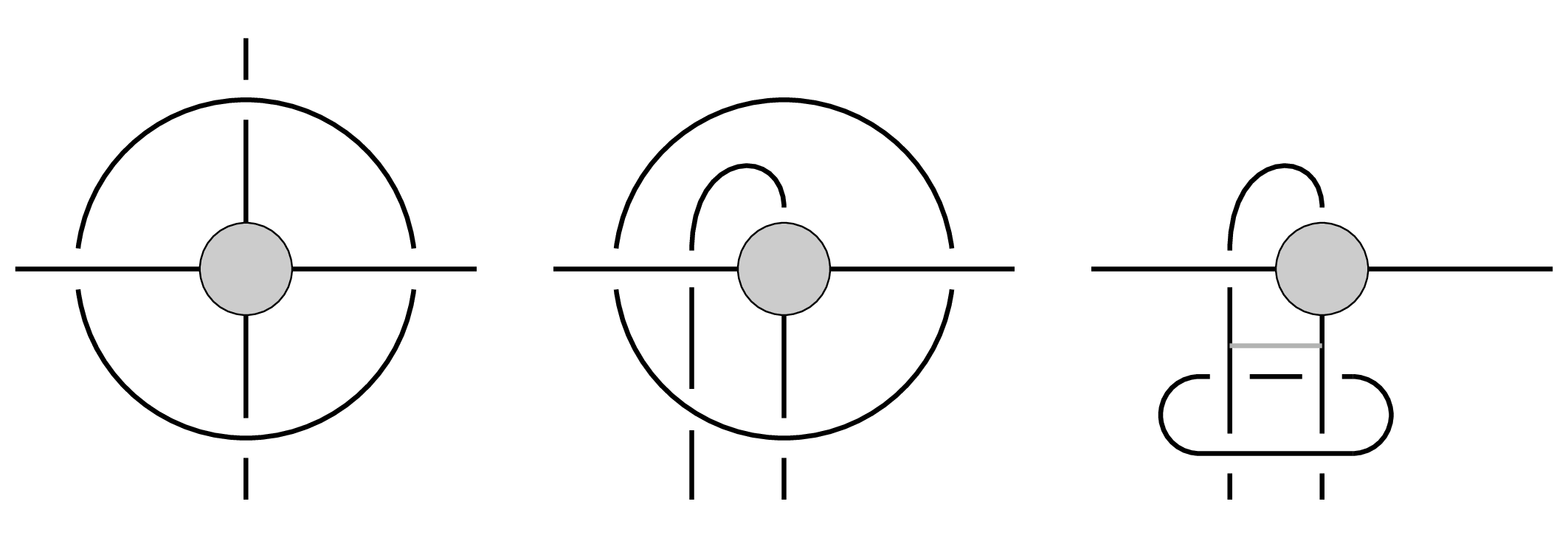}
\caption{}
\end{figure}

It is perhaps enlightening to give a heuristic, more geometric explanation of this phenomenon.
The sequence of diagrams in Fig.\ha 14 shows that the encircling simple closed curve bounds
a disk \( \Delta \) meeting \( L_i \) in \( C \) together with two transverse ``punctures''.
\( \Delta - L_i \) is thus a 3--punctured sphere in the link complement, and it follows
from Example 6.3 above that the geodesic marked \( \gamma \) in the final diagram of
Fig. 14 has associated label \( \pm\frac{1}{4} \), independent of the link \( (S^3 \bc L_i) \)
containing \( (B \bc T) \).  This imposes an extra constraint (of complex dimension 1) on the
geometric structure of the sub-tangle \( U \), making it plausible that its structure
is determined uniquely.

The boundary labels of an alternating encircled tangle are subject to further constraints,
as expressed in the next theorem.  With a little more effort the conclusion can be established
without the hypothesis that the tangle be alternating; however, our main interest here is with
alternating links.

\begin{thm}
Let \( (B \bc T) \) be an alternating encircled tangle, represented by a
standard diagram with boundary labels \( u_i \bc w_i \; (1 \leq i \leq 4) \), as in Fig. 13\,(ii).
\v

{\bf (i)}\h\h The boundary crossing labels \( w_i \) are all equal up to sign, any two being
equal if and only if the crossings to which they belong have equal sign (see Fig.\ha 15).
\v

{\bf (ii)}\h\h Opposite boundary edge labels are equal, i.e. in Fig.\ha 13\,(ii) \( u_1 = u_3 \) and
\( u_2 = u_4 \).
\end{thm}

\begin{figure}[ht!]
\labellist
\pinlabel $\mathrm{+1}$ at -110 33
\pinlabel $\mathrm{-1}$ at -110 -167
\pinlabel $\mathrm{w}$ at 42 20
\pinlabel $\mathrm{-w}$ at 305 20
\pinlabel $\mathrm{w}$ at 412 20
\pinlabel $\mathrm{w}$ at 667 20
\pinlabel $\mathrm{w}$ at 154 126
\pinlabel $\mathrm{w}$ at 524 126
\pinlabel $\mathrm{-w}$ at 196 -126
\pinlabel $\mathrm{w}$ at 555 -126
\pinlabel {\small Crossing signs} at -100 -220
\endlabellist
\centering

\includegraphics[scale=0.43]{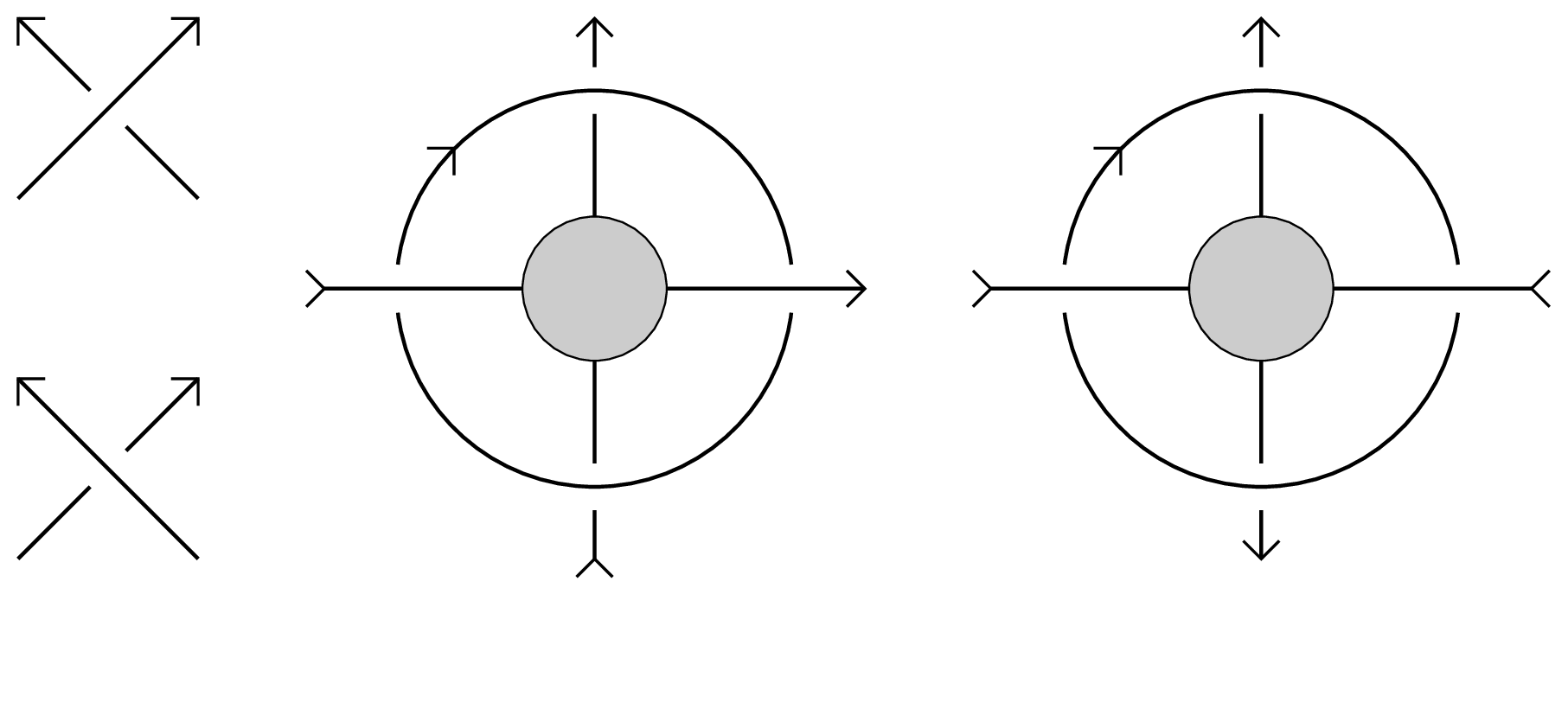}
\caption{}
\end{figure}

\begin{proof}
We first observe that conclusions (i), (ii) hold for the oriented encircled tangles
\( V_1 \bc V_2 \) illustrated in Fig.\ha 16.  This may easily be verified from the label equations
associated to the 4--sided region interior to the tangle.
For example, for \( V_1 \) we have \( u_1 = u_3 = 0 \;,\; w_1 = w_2 \;,\; w_3 = w_4 \),
together with the equations
\[ \frac{w_3}{u_2} + \frac{w_3}{u_4} = 1 \quad,\quad -\frac{w_1}{u_2} - \frac{w_1}{u_4} = 1
\quad,\quad \frac{w_3}{u_2} - \frac{w_1}{u_2} = 1 \quad,\quad \frac{w_3}{u_4} - \frac{w_1}{u_4} = 1 \h;\]
it then follows that \( u_2 = w_3 - w_1 = u_4 \) and \(\dis w_3 = -w_1 = \frac{u_2}{2} \).
If we switch the crossings of \( V_1 \), the edge labels marked \( 1 \)
both become \( -1 \), but the result of the computation is unaffected.
In the case of \( V_2 \), it is shown similarly that \( u_2 = u_4 \) and
\(\dis w_1 = w_2 = w_3 = w_4 = -\frac{u_2}{2} \).

Now suppose that \( T \) is an encircled, oriented tangle distinct from \( V_1 \bc V_2 \).  Then
we may form an alternating hyperbolic link \( L \) by summing \( T \) with a trivial tangle
(Fig.\ha 16).
We then see that \( T \) shares its encircling link component with a copy of one of
\( V_1 \bc V_2 \).  The conclusion for \( T \) then follows from the analysis of the \( V_i \)
in the previous paragraph, together with Theorem 7.3 (ii).
\end{proof}

\begin{figure}[ht!]
\labellist
\small
\pinlabel $\mathrm{u_{2}}$ at -134 37
\pinlabel $\mathrm{u_{2}}$ at 116 37
\pinlabel $\mathrm{u_{4}}$ at -63 -37
\pinlabel $\mathrm{u_{4}}$ at 187 -37
\pinlabel $\mathrm{w_{1}}$ at -124 91
\pinlabel $\mathrm{w_{1}}$ at 126 91
\pinlabel $\mathrm{w_{3}}$ at -76 -89
\pinlabel $\mathrm{w_{3}}$ at 174 -89
\pinlabel $\mathrm{1}$ at -88 27
\pinlabel $\mathrm{1}$ at 162 27
\pinlabel $\mathrm{1}$ at -121 -11
\pinlabel $\mathrm{1}$ at 129 -11
\pinlabel $\mathrm{V_{1}}$ at -97 -133
\pinlabel $\mathrm{V_{2}}$ at 153 -133
\endlabellist
\centering

\includegraphics[scale=0.39]{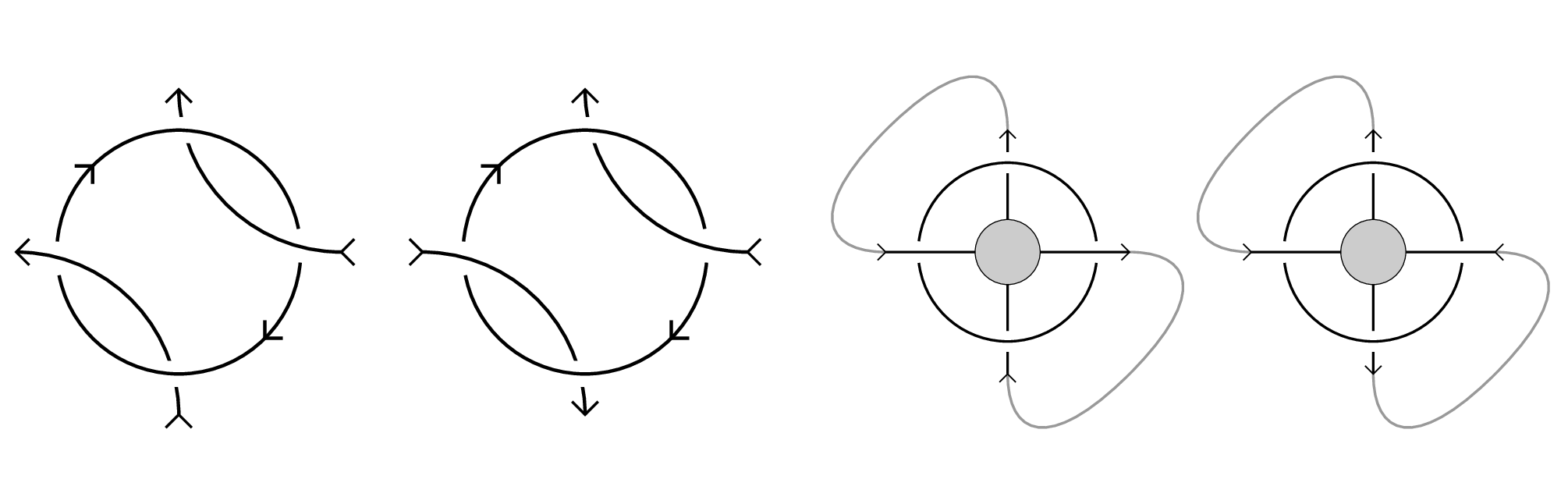}
\caption{}
\end{figure}

\end{section}

\begin{section}{Epilogue}

It had been observed by the first author that the horoball patterns of alternating links had certain ``elements
of predictability'', and that in many cases one could even reconstruct a diagram for the
link by visual inspection of the horoball pattern.  The original purpose of the method
outlined in this paper was to search for an explanation of this behaviour.  For alternating
diagrams the hyperbolicity equations relate to the diagram in a natural way, and it is perhaps
reasonable to hope that the method will lead eventually to a deeper understanding of the
geometry of prime alternating links.
\end{section}

\begin{section}{Acknowledgments}

The authors are grateful to Colin Adams, Francis Bonahon, David Futer, Darren Long and Jessica
Purcell for valuable conversations.  Particular thanks are due to the referee, whose
extremely helpful and detailed comments have played an essential part in the making
of this paper. 
\end{section}

%%%%%%%%%%%%%%%%%%%%   End of main body of article
%
%                             References
%
%   BiBTeX users uncomment the following line:
%
\bibliographystyle{gtart}

\bibliography{hyp_gt.bib}

\begin{thebibliography}{}

\bibitem[Adams, 1985]{Ad3}
Adams, C. (1985).
\newblock Thrice punctured spheres in hyperbolic 3-manifolds.
\newblock {\em Trans. Amer. Math. Soc.}, 287: 645--656.

\bibitem[Adams, 2002]{Ad1}
Adams, C. (2002).
\newblock Waist size for cusps in hyperbolic 3-manifolds.
\newblock {\em Topology}, 41 no.2: 257--270.

\bibitem[Adams, 2007]{Ad2}
Adams, C. (2007).
\newblock Noncompact fuchsian and quasi-fuchsian surfaces in hyperbolic
  3-manifolds.
\newblock {\em Algebraic and Geometric Topology}, 7:565--582.

\bibitem[Aitchison, Lumsden and Rubinstein, 1992]{AR}
Aitchison, I., Lumsden, E. and Rubinstein, H. (1992).
\newblock Cusp structures of alternating links.
\newblock {\em Invent. Math.}, 109 no.1: 473--494.

\bibitem[Bers, 1960]{Bers}
Bers, L. (1960).
\newblock Quasiconformal mappings and Teichm\"{u}ller's theorem.
\newblock {\em Analytic Functions} (R. Nevanlinna {\em et al.}, eds.),
Princeton University Press: 89--119.

\bibitem[Canary and McCullough, 2004]{CM}
Canary, R. and McCullough, D. (2004).
\newblock Homotopy equivalences of 3-manifolds and deformation theory of Kleinian groups.
\newblock {\em Memoirs Amer. Math. Soc.}, 172 no. 812.

\bibitem[Conway, 1967]{Con}
Conway, J. (1967).
\newblock An enumeration of knots and links, and some of their algebraic
  properties.
\newblock {\em Computational Problems in Abstract Algebra (Ed. Leech)}, pages
  329--358,
\newblock Pergamon Press.

\bibitem[Epstein and Penner, 1988]{EP}
Epstein, D. and Penner, R. (1988).
\newblock Euclidean decompositions of noncompact hyperbolic manifolds.
\newblock {\em J. Diff. Geom.}, 27 no.1: 67--80.

\bibitem[Francaviglia, 2004]{Franc}
Francaviglia, S. (2004).
\newblock Hyperbolic volume of representations of fundamental groups of cusped 3-manifolds.
\newblock {\em Int. Math. Res. Not.}, 9: 425--459.

\bibitem[Futer, Kalfagianni and Purcell, 2012]{FKP}
Futer, D., Kalfagianni, E. and Purcell, J. (2012).
\newblock Quasifuchsian state surfaces.
\newblock {\em Trans. Amer. Math. Soc.}, to appear.

\bibitem[Lickorish, 1981]{Lickorish}
Lickorish, W.B.R. (1981).
\newblock Prime knots and tangles.
\newblock {\em Trans. Amer. Math. Soc.}, 267 no.1: 321--332.

\bibitem[Marden, 1977]{Marden}
Marden, A. (1977).
\newblock Geometrically finite Kleinian groups.
\newblock {\em Discrete groups and automorphic functions.} (W.J. Harvey, ed.),
  Academic Press: 259--294.

\bibitem[Menasco, 1984]{Men}
Menasco, W. (1984).
\newblock Closed incompressible surfaces in alternating knot and link
  complements.
\newblock {\em Topology}, 23 no.1: 37--44.

\bibitem[Morgan, 1984]{Morgan}
Morgan, J.W. (1984).
\newblock On Thurston's uniformization theorem for three-dimensional manifolds.
\newblock {\em The Smith conjecture} (H. Bass and J.W. Morgan, eds.),
Academic Press: 37--138.

\bibitem[Menasco and Thistlethwaite, 1993]{MT}
Menasco, W. and Thistlethwaite, M. (1993).
\newblock The classification of alternating links.
\newblock {\em Annals Math.}, 138: 113--171.

\bibitem[Rolfsen, 1990]{Rolf}
Rolfsen, D. (1990).
\newblock {\em Knots and Links}.
\newblock Publish or Perish, Inc.

\bibitem[Ruberman, 1987]{Ruberman}
Ruberman, D. (1987).
\newblock Mutation and volumes of knots in $S^3$.
\newblock {\em Invent. Math.}, 90: 189--215.

\bibitem[Sakuma and Weeks, 1995]{SW}
Sakuma, M. and Weeks, J. (1995).
\newblock Examples of canonical decompositions of hyperbolic link complements.
\newblock {\em Japan J. Math. (N.S.)}, 21 no.2: 393--439.

\end{thebibliography}

\end{document}